\documentclass{amsart}

\usepackage{amssymb,latexsym,amsxtra,amscd,ifthen}

\usepackage{amsmath}

\usepackage{amsfonts}

\usepackage{verbatim}

\usepackage{dsfont}

\usepackage[capitalise]{cleveref}

\newcommand{\basis}{Gr{\"o}bner basis}
\newcommand{\A}{\widetilde{A}}

\numberwithin{equation}{section}

\theoremstyle{plain}

\newtheorem{theorem}{Theorem}[section]

\newtheorem{lemma}[theorem]{Lemma}

\newtheorem{proposition}[theorem]{Proposition}

\newtheorem{question}[theorem]{Question}

\newtheorem{corollary}[theorem]{Corollary}

\theoremstyle{definition}

\newtheorem{definition}[theorem]{Definition}

\theoremstyle{remark}

\begin{document}

\title{A classification of relation types of Ore extensions of dimension 5}
\author{Susan Elle}

\keywords{Artin-Schelter regular, Ore, enveloping algevra, graded Lie algebra, Hilbert series, noncommutative projective geometry, noetherian graded ring, Grobner basis}

\subjclass[2010]{Primary: 16S36, 16S38; Secondary: 14A22, 16U20, 16W50, 16-04}

\thanks{The author was partially supported by NSF grant DMS-0900981.}
\address{Susan Elle\\
	Department of Mathematics\\
	University of California, San Diego\\
	La Jolla, CA 92093-0112\\
	USA}

\email{selle@ucsd.edu}

\begin{abstract}
	In order to study AS-regular algebras of dimension 5, we consider dimension 5 graded iterated Ore extensions generated in degree one.  We classify the possible degrees of relations and structure of the free resolution for extensions with 3 and 4 generators.  We show that every known type of algebra of dimension 5 can be realized by an Ore extension and we consider which of these types cannot be realized by an enveloping algebra.
\end{abstract}
\maketitle
\section{Introduction}
The study of Artin-Schelter (AS) regular algebras was introduced by Artin and Schelter in 1987 \cite{MR917738}.  AS-regular algebras correspond to noncommutative homogeneous coordinate rings of weighted projective spaces and so their classification is of interest in the field of noncommutative algebraic geometry.  Our goal in this paper is to classify possible types of graded iterated Ore extensions of dimension 5 which are generated in degree one.  We present an interesting example of an Ore extension with 2 degree one generators with the property that it has a Hilbert series which cannot be realized by any enveloping algebra. We also list all possible types of dimension 5 iterated Ore extensions with 3 and 4 degree one generators and consider which of these cannot be realized by the enveloping algebra of any $\mathds{N}$-graded Lie algebra.

An Artin-Schelter regular algebra over K that is generated in degree one has finite presentation $\displaystyle A=\frac{K\langle x_1,\cdots,x_b\rangle}{I}$, and its trivial module, $K$, has minimal free resolution:
\begin{multline} \label{Resolution}
0\rightarrow A(-l)\rightarrow A(-l+1)^b\rightarrow \bigoplus \limits_{i=1}^{n} A(-l+a_i)\rightarrow\cdots \\
\cdots\rightarrow \bigoplus \limits_{i=1}^{n}  A(-a_i)\rightarrow A(-1)^b\rightarrow A\rightarrow K\rightarrow 0,
\end{multline}
where any minimal generating set of $I$ has $n$ elements with the degree of the $i$th element equal to $a_i$. (For a proof see, for example, \cite[Lemma 2.11 together with the discussion following Definition 2.1]{R}.)  This resolution is often described via \emph{graded Betti numbers} where $\beta_{i,j}$ is equal to the number of copies of $A(-j)$ appearing in the $ith$ step of the resolution.  The \emph{Hilbert series} of A is $h_A(t)=\sum\limits_{n\in \mathds{N}} (\dim_KA_n)t^n$ where $A_n$ is the $nth$ graded piece of $A$.  The Hilbert series can be computed directly from the free resolution: $\displaystyle h_A(t)=\frac{1}{q(t)}$ where $\displaystyle q(t)=\sum\limits_{i,j} (-1)^i\beta_{i,j}t^j$. (See, for example, \cite[2.6]{R}.)

 Thus, there are many invariants that we can use to discuss the possible \emph{classification} of types of AS-regular algebras.  Most generally, we can use their Hilbert series, although there are algebras with fundamentally different structures that share the same series.  More refined, we can use their \emph{relation type} (the number and degree of the relations in the minimal generating set of $I$, which we will denote by $(a_1,\cdots, a_n)$ where $a_1\leq \cdots\leq a_n$).  More refined still, we can refer to the \emph{resolution type}, or the set of graded Betti numbers of A. The most concrete option for classifying AS-regular algebras, and one beyond the scope of this paper, would be to list the possbile \emph{families of relations} for the algebras by explicitly writing the possible coefficients of the relations.  For example, an AS-regular algebra of dimension 2 which is generated in degree one is isomorphic to $\displaystyle\frac{K\langle x_1,x_2\rangle}{\langle r \rangle}$ where $r=x_2x_1-qx_1x_2$, $0\neq q$ (in which case the algebra is called the \emph{quantum plane}) or $r=x_2x_1-x_1x_2-x_1^2$ (and the algebra is called the \emph{Jordan plane}).
 
 Substantial progress has been made on the classification of AS-regular algebras generated in degree one of dimension 3 and 4 \cite{MR1086882, MR1128218, MR2309153,RZ, MR2452318, MR2529094}.  In dimension 3, the possible families of relations are known.  In dimension 4, the possible resolution types are known provided that the algebra is a domain, and under mild assumptions, the possible families of relations are known in the case where the algebras are also assumed to be \emph{$\mathds{Z} ^2$-graded}, i.e. $\displaystyle A=\frac{K\langle x_1,\cdots,x_b\rangle}{I}$, $deg(x_i)\in\{(1,0),(0,1)\}$ for all $i$, and $I$ is homogeneous in $\mathds{Z}\times \mathds{Z}$.  In dimension less than 5, it is known that each Hilbert series has a unique resolution type, every resolution type can be realized by the universal enveloping algebra of an $\mathds{N}$-graded Lie algebra, and every resolution type can be realized by an algebra which is $\mathds{Z} ^2$-graded.
 
 More recently, Floystad and Vatne studied dimension 5 AS-regular algebras generated in degree one with 2 generators.  They found an example of an AS-regular algebra for which there is no enveloping algebra with the same Hilbert series \cite[Section 4]{FV} and, under mild assumptions, provided a short list of all possible relation types, leaving open the question of whether there exist AS-regular algebras with relation type (4,4,4,5) or (4,4,4,5,5) [Section 5].  Wang and Wu used $A_\infty$ techniques to further classify dimension 5 algebras with 2 generators.  They found several families of algebras with relation type (4,4,4,5,5) and proved that there is an Ore extension of this type (although there is no enveloping algebra with this relation type) \cite[Section 5.2]{WW}.  Zhou and Lu further classified the algebras under the additional assumption of a $\mathds{Z} ^2$-grading and listed all possible families of relations in each case.  They found there is no $\mathds{Z} ^2$-graded algebra of type (4,4,4,5) \cite[Proposition 5.3]{ZL}, although it remains an open question whether a regular algebra with this relation type exists.\\
 
 To classify Ore extensions we note that, by the symmetry of the free resolution, the relation type uniquely determines the resolution type in dimension 5, although this is probably false in higher dimensions.  After introducing some preliminary definitions and a convenient presentation for Ore extensions in Section 1 and listing the possible degrees of variables of a 5-dimensional Ore extension in Section 2, we prove the following in Section 3.
 
 \begin{theorem}	
 	[\Cref{Ore11235}] There is an Ore extension with 2 degree one generators that has Hilbert series $\displaystyle h_A(t)= \frac{1}{(1-t)^2(1-t^2)(1-t^3)(1-t^5)}$ and relation type (3,4,7).
 \end{theorem}
 Together with other results in the field, this means that every known type of AS-regular algebra can be realized by an Ore extension.  It remains an open question whether there is an AS-regular algebra of relation type (4,4,4,5) and an Ore extension of the same type, although our preliminary calculations suggest that the latter is unlikely.  In Section 4 we prove that:
 
 \begin{theorem} [\Cref{No(3333)}, \Cref{No(333)}, \Cref{EnvelopingExamples}]
 	The relation types for a dimension 5 Ore extension with 4 degree one generators are (2,2,2,2,2), (2,2,2,2,2,3), and (2,2,2,2,2,3,3).  There is no enveloping algebra with type \linebreak(2,2,2,2,2,3) although there are enveloping algebras for the other relation types.
 \end{theorem} 
 In Section 5, we prove:
 \begin{theorem} [\Cref{(11122)}, \Cref{Env(11122)}, \Cref{(11123)}, \Cref{Env(11123)}]
 	The relation types for a dimension 5 Ore extension with 3 degree one generators are (2,3,3,3,3,3), (2,2,3), and (2,2,3,4).  There is no enveloping algebra with type (2,2,3) although there is an enveloping algebra for types (2,2,3,4) and (2,3,3,3,3,3).
 \end{theorem}
To complete the classification of Ore extensions of dimension 5, it is then natural to ask:
 \begin{question}
 	Is there a dimension 5 iterated Ore extension with 2 degree one generators and relation type (4,4,4,5)?
 \end{question}

Since our classification of Ore extensions was motivated by a desire to classify AS-regular algebras in general, we also would like to know:
\begin{question}
	Is there an AS-regular algebra of dimension 5 with one of the Hilbert series considered in this paper that has a different relation type than that of an Ore extension? 
\end{question}

\begin{question}
	Is there an AS-regular algebra of dimension 5 with a different Hilbert series than those considered in this paper? 
\end{question}

Both of these would help to answer our larger underlying question:
\begin{question}
	Can every type of AS-regular algebra be realized by an Ore extension?
\end{question}
%\[0\rightarrow \bigoplus\limits_{j\in\mathds{N}}A(-j)^{\beta_{d,j}} \rightarrow \cdots \rightarrow \bigoplus \limits_{j \in \mathds{N}} A(-j)^{\beta_{2,j}} \rightarrow \bigoplus \limits_{j \in \mathds{N}} A(-j)^{\beta_{1,j}} \rightarrow A \rightarrow K \rightarrow 0\]

\section{Preliminaries}
In this section we review basic definitions for AS-regular algebras and iterated Ore extensions and find a convenient presentation for the latter.  The interested reader may find \cite[Chapters 1 and 2]{MR2080008} a useful reference for some of the claims made about Ore extensions.

\begin{definition}\label{DefOre}
Let $R$ be a ring. An \emph{Ore extension $ R[x,\sigma,\delta]$} is a ring with elements of the form $f(x)=\sum\limits_{i=0}^n a_ix^i, a_i\in R$ and multiplication satisfying $xr=\sigma(r) x+\delta(r)$ for all $r \in R$, where $\sigma$ is an endomorphism of $R$ and $\delta$ is a \emph{$\sigma$-derivation} of $R$,
 i.e. $\delta(r_{1}r_{2})=\sigma (r_{1})\delta (r_{2})+\delta (r_{1})r_{2}$ for all $r_1,r_2\in R$.

An \emph{iterated Ore extension} $ R[x_1,\sigma_1,\delta_1][x_{2},\sigma _{2},\delta _{2}]\dotsm [x_{n},\sigma _{n},\delta _{n}]$ is an Ore extension where for all $j\geq1$, $\sigma_j$ and $\delta_j$ are a ring endomorphism and a $\sigma_j$-derivation of \emph{$R_{(j-1)}$ }$:= R[x_1,\sigma_1,\delta_1][x_{2},\sigma _{2},\delta _{2}]\dotsm [x_{j-1},\sigma _{j-1},\delta _{j-1}]$, respectively.  Elements in this extension have the form $\sum\limits_{i=0}^n a_ix_1^{i_1}\cdots x_n^{i_n}, \; a_i\in R.$

We are most interested in specific types of iterated Ore extensions and will modify our definition as additional desirable properties are introduced.
\end{definition}

%%%%%%%%%%%%%%%%%%%%%%%%%%%%%%%%%%%%%%%%%%%%%%%
%BEGIN_FOLD
\iffalse
\begin{definition}
A \emph{degree function} on a ring $R$ is a function $d:R\rightarrow \mathds{N}\cup\{ -\infty \}$ such that for all $f$ in $R$
\begin{enumerate}
\item $d(f) \geq 0$ for $f \neq 0$, $d(0)=-\infty$,
\item $ d(f-g) \leq max\{d(f),d(g)\} $, and
\item $d(fg)=d(f)+d(g)$.
\end{enumerate}
\end{definition}
\fi
%END_FOLD	

We now find it convenient to review some definitions and notation.  A more thorough introduction to the material, as well as a proof of \Cref{Bergman}, can be found in Bergman's paper \cite[Section 1]{B}.

 Suppose $R$ is an associative algebra with unity over a field $K$ and we have a presentation of $R$ by a family $X$ of generators and a family $S$ of relations. In practice, we care about the ideal generated by $S$, call it $I$. We have $R\cong \displaystyle \frac{K\langle X\rangle}{I}$ where $K\langle X\rangle$ is the free associative K-algebra on $\langle X\rangle$ and $\langle X\rangle$ is the free semigroup with $1$ on $X$ .  A subset, $B\subseteq S$, is a \emph{minimal generating set} for $I$ if no proper subset of $B$ generates $I$.  Fix a total ordering on $\langle X\rangle$ with the property that if $w<v$ then $uw<uv$ and $wu<vu$ for all $u\in \langle X\rangle$. Such an ordering will be called a \emph{semigroup total ordering}. Every relation $\sigma \in S$ can be written in the form $W_\sigma=f_\sigma$ where $W_\sigma$ is a monomial and is larger than any of the monomials in $f_\sigma$.  We call $W_\sigma$ the \emph{leading term} (denoted LT) of the relation $\sigma$.  We can assume that the leading term is always monic since $K$ is a field.  We can also take S such that all leading terms are distinct (since otherwise we could subtract a scalar multiple of one relation from another to get two relations with different leading terms which generate the same ideal).

A word $w$ is \emph{irreducible} under $S$  if it does not contain any $W_\sigma$ as a subword.  Otherwise, $w$ contains some $W_\sigma$, say $w=uW_\sigma v$ and we consider the $K$-linear \emph{reduction} map $r_{uW_\sigma v}:K\langle X\rangle \rightarrow K\langle X \rangle$ which sends $uW_\sigma v$ to $uf_\sigma v$ and fixes all other elements of $\langle X \rangle$. A finite sequence of reductions $r_1\cdots r_n$ ($\displaystyle r_i=r_{u_i W_{\sigma_i}v_i}$) is \emph{final} on $w$ if $r_1\cdots r_n(w)$ is irreducible. The word $w$ is \emph{reduction unique} if its images under all final sequences of reductions are the same.

A $5$-tuple $(\sigma,\tau, u,v,w)$ with $\sigma,\tau\in S$ and $u,v,w \in {\langle X \rangle}$ is an \emph{overlap ambiguity} if $u, v, w\neq$1, $W_\sigma =uv$, and $W_\tau =vw$ and an \emph{inclusuion ambiguity} if $\sigma\neq\tau$,
$W_\sigma=v$, and $W_\tau=uvw$. An ambiguity is \emph{resolvable} if there exist compositions of reductions \linebreak $s$ and $s'$ such that $s(r_{W_\sigma w}(uvw))=s'(r_{uW_\tau}(uvw))$ (in the case of an overlap ambiguity) or $s(r_{uW_\sigma w}(uvw))=s'(r_{W_\tau}(uvw))$ (in the case of an inclusion ambiguity). A set $S$ of relations satisfies the \emph{diamond condition} if all reduction ambiguities (overlap and inclusion) are resolvable, in which case we say that $S$ is a \emph{Gr{\"o}bner basis} of $I$.  A \basis $ $ is \emph{reduced} if all leading terms are monic and no element in the basis has a monomial which contains the leading term of any other element in the basis.  Throughout the rest of this paper, any reference to a \basis $ $ will mean a reduced \basis.

\begin{theorem} \cite[Theorem 1.2]{B}\label{Bergman}
Let $\leq$ be a semigroup total ordering having the descending chain condition and let $S$ be a set of relations where the leading term of each relation is monic and distinct from the leading term of any other relation.  Then the following conditions are equivalent:
\begin{enumerate}
\item S satisfies the diamond condition;
\item All elements of $K\langle X\rangle$ are reduction unique under $S$;
\item A set of representatives for the elements of the algebra $\displaystyle R=\frac{K\langle X \rangle}{I}$ determined by the generators X and the ideal $I$ generated by the relations S is given by the $K$-submodule $K\langle X \rangle_{irr}$ spanned by the $S$-irreducible monomials of $\langle X \rangle$.
\end{enumerate}
\end{theorem}

Continuing with definitions, we shall say that an algebra $R$ over $K$ is \emph{$\mathds{N}$-graded} if $R=\bigoplus \limits_{n=0}^\infty R_n$ as $K$-spaces and $R_nR_m\subseteq R_{n+m}$ for all $n$ and $m$.  $R$ is \emph{connected} if $R_0=K$.  A \emph {graded iterated Ore extension} is an iterated Ore extension as in \Cref{DefOre} with variables $(x_1,x_2,\cdots,x_n)$ of degrees $(deg(x_1),deg(x_2),\cdots,deg(x_n))$, $deg(x_i)\geq 1$, where $\sigma_j(x_i)$ is homogeneous of degree $deg(x_i)$ and $\delta_j(x_i)$ is homogeneous of degree $deg(x_j)+deg(x_i)$ for all $n\geq j>i\geq1$.  In particular, such a ring is $\mathds{N}$-graded.  We refer to $(deg(x_{i_1}),\cdots,deg(x_{i_n}))$ as the \emph{degree type} of an Ore extension and require that the degrees be listed in ascending order so that the expression is unique.  For the rest of this paper, any mention of an \emph{Ore extension} will refer to a \emph{graded iterated Ore extension} over $K$ unless otherwise stated. 

\emph{Graded lexicographic order} is a total order on $\langle X\rangle$ where $w_1>w_2$ if $deg(w_1)>deg(w_2)$ or if $deg(w_1)=deg(w_2)$ and $w_1=a_1a_2\cdots a_j$ comes before $w_2=b_1b_2\cdots b_k$ in the lexicographic order.  For the rest of this paper, we will use graded lexicographic order.  This has the descending chain condition.  It is sometimes convenient to consider the case when the (lexicographic) order is taken to be
\[ x_{n}>x_{n-1}>\cdots >x_{1},\]
and we will assume that this is the order on the variables for the rest of this section.

A finitely generated \emph{$\mathds{N}$-graded Lie algebra} is a Lie algebra with generators \linebreak $\{x_1,\cdots, x_n\}$ assigned positive degrees such that the Lie bracket preserves degree, i.e. each monomial in $[x_j,x_i]$ has degree equal to $deg(x_j)+deg(x_i)$.  We note without proof that the universal enveloping algebra of an $\mathds{N}$-graded finite dimensional Lie algebra $L$ with $K$-basis $\{x_1,\cdots ,x_n\}$ can be taken to have $deg(x_1)\geq\cdots\geq deg(x_n)$ and lexicographic order $x_n\geq\cdots\geq x_1$ and will then have presentation
 $\displaystyle U=\frac {K\langle x_{1}\cdots x_{n}\rangle } {\langle \{r_{ji} \} \rangle }$ where for each $j>i$, there is a unique homogeneous relation $r_{ji}$ given by 
 \[r_{ji}:x_{j}x_{i}=x_i x_j+ \hspace{-6 mm} \sum_{\substack{k \, | \, deg(x_{k})=\\ deg(x_j )+deg(x_i )}}  \hspace{-6 mm} a_{ji}^k x_{k}, \; a_{ji}^k \in K \]
 and where the relations satisfy the diamond condition.  This is a consequence of the  Poincar{\'e}-Birkhoff-Witt theorem and is proven in Bergman's paper \cite[Theorem 3.1]{B}.  The interested reader may also wish to refer to Humphrey's book on Lie algebras \cite[Chapter 1 and Chapter 17]{MR0323842} for the relevant background information.  This statement has a converse which also holds.
\begin{theorem}\label{PresentationEnv}
	If $U$ has presentation $\displaystyle U=\frac {K\langle x_{1}\cdots x_{n}\rangle } {\langle \{r_{ji} \} \rangle }$ where for each $j>i$, there is a unique homogeneous relation $r_{ji}$ given by 
	\[r_{ji}:x_{j}x_{i}=x_i x_j+ \hspace{-6 mm} \sum_{\substack{k \, | \, deg(x_{k})=\\ deg(x_j )+deg(x_i )}}  \hspace{-6 mm} a_{ji}^k x_{k}, \; a_{ji}^k \in K \]
	and where the relations satisfy the diamond condition, then $U$ is the universal enveloping algebra of a graded Lie algebra.
	\begin{proof}
		Suppose an algebra $U$ has the presentation described.  Define $L$ to be generated as a $K$-vector space by $\{x_1,\cdots , x_n \}$ and define a multiplication on the generators of $L$ by
		\begin{displaymath}
		[x_j,x_i] = \left\{
		\begin{array}{lr}
		\displaystyle \sum_k a_{ji}^k x_{k}  &j>i,\\
		\displaystyle \sum_k -a_{ji}^k x_{k}  &j<i,\\
		0 & j=i.
		\end{array}
		\right.
		\end{displaymath}
		This multiplication  can be extended bilinearly to general elements in $L$.  We claim that $L$ is the desired graded Lie algebra.  The multiplication satisfies bilinearity and has the alternating property by construction.  That the Jacobi identity is satisfied is equivalent to the fact that all ambiguities in $U$ resolve by the proof of theorem 3.1 in Bergman's paper \cite{B}.  Finally, since this multiplication is degree preserving, $L$ is a graded Lie algebra with enveloping algebra $U$.	
	\end{proof}
\end{theorem}

Since both enveloping algebras and Ore extensions have the same basis as a weighted commutative polynomial ring with the same variables and degrees, the Hilbert series of these algebras is known: $\displaystyle h(t)=\frac{1}{\prod\limits_{i=1}^n (1-t^{deg(x_i)})}$.

What follows are slightly modified versions of the theorems in Cohn's book Algebra Vol. 2 \cite[Chapter 12, Theorem 1]{C} and the proofs follow from the original proof.  If $R$ is a domain, $\sigma$ a degree preserving endomorphism, and $\delta$ a degree preserving $\sigma $-derivation, then there exists an Ore extension $P=R[x,\sigma,\delta]$.  Conversely we have:

\begin{theorem} \label{Cohn2}
	Let $R$ be a non-trivial graded ring and let $P$ be a graded ring containing $R$ with element $x\in P$ such that elements of $P$ can be written uniquely in the form $f(x)=\sum\limits_{i=0}^n a_ix^i$, $a_i\in R$, and satisfy a homogeneous relation $xa=\sigma (a)x+\delta (a)$ for some $\sigma (a), \delta(a) \in R$.  Then $\sigma$ is an endomorphism, $\delta$ a $\sigma$-derivation, and $P\cong R[x,\sigma,\delta]$ is an Ore extension.  
	
	%BEGIN_FOLD
	%%%%%%%%%%%%%%%%%%%%%%%%%%%%%%%%%%%%%%%%%%%%%%%%%%%%
	\iffalse
	% Note: changes required since \sigma need not be injective
	
	\begin{proof}
		We can compute
		\begin{align*}
		x(a+b) &= \sigma (a+b)x+\delta (a+b)\\
		x(a+b) &= xa+xb\\
		&= \sigma (a)x+ \delta(a) + \sigma(b)x+ \delta(b)\\
		&= (\sigma (a)+\sigma (b))x+(\delta (a)+\delta (b), \text{ while}\\
		x(ab) &= \sigma (ab)x+\delta (ab)\\
		x(ab) &=(\sigma (a)x+\delta (a))b\\
		&= (\sigma (a)\sigma (b))x+(\sigma (a)\delta (b)+ \delta (a)b).
		\end{align*}
		By the uniqueness of the form, $\sigma$ must be an $R$ homomorphism and $\delta$ a $\sigma $-derivation.
		
		If $P$ also has the given degree function, $R$ must be entire.  For if $ab=0$ and $deg(fg)=deg(f)+deg(g)$, then $(xa)(bx)=0$ forces either $a=0$ or $b=0$.  Similarly, $\sigma$ must be injective.  For if $\sigma(a)=0$ then $deg(xa)= deg(\sigma(a)x +\delta(a))$ forces $a=0$.  Thus, if $P$ has the given degree function, $P=R[x;\sigma,\delta]$ is an Ore extension.
		
	\end{proof}
	\fi
	%END_FOLD
	%%%%%%%%%%%%%%%%%%%%%%%%%%%%%%%%%%%%%%%%%%%%%%%%%%%%%%%%%%%%%%
	
\end{theorem}

We now consider a presentation for (graded iterated) Ore extensions over a field $K$.  Since we want these to be $K$-algebras, we are interested in the case where $K$ is central, which means that $\sigma_1$ is the identity and $\delta_1$ is the zero mapping.  We note that the following two theorems would also hold for ungraded iterated Ore extensions with the term ``homogeneous" removed from the proofs, but these are of lesser interest to us.  Recall our notation established in \Cref{DefOre}: $K_{(j)}:=K[x_1]\cdots[x_j,\sigma_j,\delta_j]$.

\begin{theorem} \label{Presentation 1}
If $K$ is a field and $ P= K[x_{1}][x_{2},\sigma _{2},\delta _{2}]\dotsm [x_{n},\sigma _{n},\delta _{n}]$ is a (graded iterated) Ore extension,  then P has presentation
\[ P\cong \frac {K\langle x_{1}\cdots x_{n}\rangle } { \langle \{r_{ji} \} \rangle }
\]
where  for each $j>i$, there is a unique homogeneous relation $r_{ji}$, given by 
\[r_{ji}:x_{j}x_{i}=\sigma_{j}(x_{i})x_{j}+\delta _{j} (x_{i}),\: \sigma _{j} (x_{i})\text{ and } \delta _{j}(x_{i}) \in K_{(j-1)},\] and these relations satisfy the diamond condition.
 
\begin{proof}
Let $K$ be a field and suppose $P$ is an iterated Ore extension:\linebreak $P=K[x_{1}][x_{2},\sigma _{2},\delta _{2}]\dotsm [x_{n},\sigma _{n},\delta _n]$.  Given $j>i$, $x_{j}x_{i}=\sigma _{j} (x_{i})x_{j}+\delta _{j} (x_{i})$ where $\sigma _{j}$ is an endomorphism of $K_{(j-1)}$ and $\delta _ {j}$ is a $\sigma $-derivation of $K_{(j-1)}$ and every monomial in the equation has degree equal to $deg(x_j)+deg(x_i)$. Since these relations allow any element in $P$ to be written as a linear combination of terms of the form $k x_{1}^{e_{1}}\cdots x_{n}^{e_{n}}$, the leading term of any additional relation would be of this form, which would contradict the fact that $\{ x_{1}^{e_{1}}\cdots x_{n}^{e_{n}}\}$ is a K-basis for the Ore extension $P$.  Thus, there cannot be any additional relations so each $r_{ji}$ is unique, all reduction ambiguities must resolve, and the diamond condition is satisfied.
\end{proof}
\end{theorem}
We can also prove a converse:
\begin{theorem} \label{Presentation2}
If $K$ is a field and $\displaystyle P=\frac {K\langle x_{1}\cdots x_{n}\rangle } {\langle \{r_{ji} \} \rangle }$ where for each $j>i$, there is a unique homogeneous relation $r_{ji}$ given by 
$x_{j}x_{i}=\sigma_{j}(x_{i})x_{j}+\delta _{j} (x_{i})$, $\sigma _{j} (x_{i}) \text{ and } \delta _{j}(x_{i}) \in \displaystyle\frac{K\langle x_1,\cdots,x_{j-1}\rangle}{\langle \{ r_{ji} \}\rangle}$, and these relations satisfy the diamond condition, then $P\cong K[x_{1}][x_{2},\sigma _{2},\delta _{2}]\dotsm [x_{n},\sigma _{n},\delta _{n}]$ is a (graded iterated) Ore extension.

\begin{proof}

Assume that we have the unique relations $\{ r_{ji}\}$ satisfying the diamond condition and for the purpose of induction assume that 
\[
R=\frac {K\langle x_{1}\cdots x_{m-1}\rangle} {\langle \{r_{ji} \} \rangle }\cong K[x_{1}][x_{2},\sigma _{2},\delta _{2}]\dotsm [x_{m-1},\sigma _{m-1},\delta _{m-1}]
\]
 is an Ore extension.  Then any monomial in $\displaystyle \frac {K\langle x_{1}\cdots x_{m}\rangle} {\langle \{r_{ji} \} \rangle }$ has the form \linebreak  $x_{m}^{c_{1}}s_{1}x_{m}^{c_{2}}s_{2}\cdots x_{m}^{c_{q}}s_{q}, $ $ s_i \in R$ where by induction each $s_{i}$ can be taken to have the form  $ k x_{1}^{f_{1}}\cdots x_{m-1}^{f_{m-1}}$ since  $\{ x_{1}^{e_{1}}\cdots x_{m-1}^{e_{m-1}}\}$ is a basis for $R$. By repeated application of the relations $x_{m}x_{i}=\sigma _{m} (x_{i})x_{m}+\delta _{m} (x_{i})$, the monomial can be written in the form $s' _{a} x_{m}^{a}+ \cdots + s'_{1}x_{m}+s' _{0}$,  $s'_{i}\in R$.  This representation is unique since the $\{r_{ji}\}$ satisfy the diamond condition by assumption.  Thus $\{x_1^{e_1},\cdots,x_m^{e_m}\}$ is a $K$-basis by \cref{Bergman}. Since the relations were also chosen to be homogeneous, $R[x_m,\sigma_m,\delta_m]$ is an Ore extension by \cref{Cohn2} and by induction, $ P\cong  K[x_{1}][x_{2},\sigma _{2},\delta _{2}]\dotsm [x_{n},\sigma _{n},\delta _{n}]$ is an Ore extension.
\end{proof}
\end{theorem}

Since Ore extensions for which $\sigma_j$ is an automorphism for all $j\geq1$ have especially nice properties, we also find the following result helpful:

\begin{theorem} In a (graded iterated) Ore extension
	$K[x_{1},\sigma_1,\delta_1]\dotsm [x_{n},\sigma _{n},\delta _{n}]$, for any $j\geq1$, if $\sigma_j$ is injective then it is an automorphism of $K_{(j-1)}$.
	
	\begin{proof} Let $\displaystyle K_{(j-1)}^i$ denote the $i$th graded piece of $K_{(j-1)}$, i.e. the set of all degree $i$ homogeneous polynomials in $K_{(j-1)}$.  $K_{(j-1)}^i$ has finite $K$-basis $\{x_1^{f_1} x_2^{f_2}\cdots x_{j-1}^{f_{j-1}} \; |  \linebreak f_1 deg(x_1)+f_2 deg(x_2)+\cdots+f_{j-1} deg(x_{j-1})=i\}$.  By the definition of a graded iterated Ore extension, we know that $\sigma_j$ preserves degree on the generators and hence on all of $K_{(j-1)}$.  Any injective map from a finite-dimensional vector space to itself must also be surjective by the rank-nullity theorem, so for any $i$ and $j$,  $ \sigma_j|_{K_{(j-1)}^i}:K_{(j-1)}^i \rightarrow K_{j-1}^i$ is bijective, and thus $\sigma_j$ is an automorphism of $K_{(j-1)}$ for all $j$.
	\end{proof}
	
\end{theorem}

Recall that we are interested in the study of Ore extensions because they provide examples of AS-regular algebras.  We provide the definition of such algebras here, but readers looking for a more detailed introduction to the material may wish to refer to the notes found in \cite[Lecture 1 and Lecture 2]{R}.

\begin{definition}
	A connected graded algebra $A=\bigoplus\limits_{i=0}A_i$ is \emph{Artin-Schelter} regular of dimension d if
	\begin{enumerate}
		\item A has finite global dimension $d$;\\
		\item A has finite Gelfand-Kirillov dimension;\\
		\item A is AS-Gorenstein, i.e.\\
		\[\text{Ext}_A^i(K,A) =\left\{
		\begin{array}{lr}
			0 &  i\neq d\\
			K(l) &  i=d\\
		\end{array}
		\right.\]
		where $K(l)$ is a shifted copy of $K$ satisfying $K(l)_n = K_{l+n}$.
	\end{enumerate}
\end{definition}

It is a fact that the universal enveloping algebra of a graded Lie algebra is Artin-Schelter regular \cite[Theorem 2.1]{FV}.  From the presentations provided, it is also clear that any such enveloping algebra is a specific example of an Ore extension.

More generally, it is also known that Ore extensions where $\sigma_j$ is an automorphism for all $j\geq1$ are Artin-Schelter regular (see \cite[Proposition 2]{MR1127037}).

Motivated by the study of AS-regular algebras, our goal in the rest of this paper is to classify the possible relation and resolution types of all dimension 5 ``Ore extensions," by which we mean ``graded iterated Ore extensions with injective (and thus bijective) $\sigma_j's$, generated in degree one."  For the sake of brevity, we will wish to have any easy way to refer to such algebras.

\begin{definition}
	An \emph{AS-Ore} extension is a graded iterated Ore extension with $\sigma_j$ injective for every $j\geq 1$ and which is generated in degree one as a $K$-algebra.
\end{definition}

We note that this definition is in no way standard (and in particular, there are AS-regular algebras that are not generated in degree one and which we have chosen not to study at this time).  We also note that the enveloping algebra of an $\mathds{N}$-graded Lie algebra which is generated in degree one is also an AS-Ore extension.

\section{A classification of possible degree types of AS-Ore extensions}

Our goal in this section is to list the 7 possible degree types for an Ore extension generated in degree one with 5 variables.  We note that when considering fully general Ore extensions, we may either order the variables by descending degree ($deg(x_1)\leq\cdots\leq deg(x_n)$) at the expense of fully controlling the lexicographic order in the Ore extension $K[x_{i_1}]\cdots[x_{i_n},\sigma_{i_n},\delta_{i_n}]$, or we may assume that $x_5>\cdots>x_1$ in the lexicographic order at the expense of controlling the degrees of these variables.  We transfer freely between these two conventions depending on which is more convenient in each situation and the convention we use does not affect the validity of any theorems we prove for general extensions. 
%
\iffalse
Here we must consider fully general Ore extensions $ K[x_{i_1}][x_{i_2},\sigma_{i_2},\delta_{i_2}]\cdots [x_{i_n},\sigma_{i_n}.\delta_{i_n}]$ with $x_{i_1}<x_{i_2}<\cdots <x_{i_n}$ in the lexicographical ordering.
\fi
%The result relies heavily on the following lemma.

\begin{lemma} \label{ijk}  If $A=  K[x_{1}][x_{2},\sigma_{2},\delta_{2}]\cdots [x_{n},\sigma_{n},\delta_{n}]$ is an AS-Ore extension and $deg(x_k)\neq1$ then there exist $i$ and $j$ with $deg(x_i)+deg(x_j)=deg(x_k)$.

\begin{proof}
Assume $deg(x_k)>1$.  If the algebra is generated in degree one, $x_k=f(\hat{X})$ where $\hat{X}=\{x\in X|deg(x)=1\}$, $X=\{x_1,\cdots, x_n\}$, and $f(\hat{X})$ is a (noncommutative) polynomial in variables from  $\hat{X}$.  This gives the relation $0=x_k-f(\hat{X})$.  By \Cref{Presentation 1}, $\displaystyle A\cong \frac{K\langle x_1,\cdots x_n\rangle}{I}$ where $I=\langle \{r_{ji}\}_{j>i} \rangle$, so any relation is generated by the $\{r_{ji}\}$ and we have that, in the free algebra, $x_k-f(\hat{X})=\sum\limits_{i,j} p_{ji}(X)r_{ji}q_{ji}(X)$ where $p_{ji}(X)$ and $q_{ji}(X)$ are (noncommutative) polynomials in $X$ and the $r_{ji}$ are the generators of $I$ and hence have degree greater than zero.

Equating polynomials, we find that the monomial $x_k$ must appear in the right side of this equation, so there exist fixed $i$ and $j$ and monomials $m_p$, $m_r$, and $m_q$ of $p_{ji}$, $r_{ji}$, and $q_{ji}$ with $x_k=m_pm_rm_q$. Since $deg(m_r)>0$, we get that $m_p$ and $m_q$ must be scalars and $ax_k$ is a monomial of $r_{ji}$ where $0\neq a\in K$. So $r_{ji}$ is a relation with leading term $x_jx_i$ and has now been shown to have a scalar multiple of $x_k$ as a term.  Since $r_{ji}$ is also homogeneous, this means that $deg(x_i)+deg(x_j)=deg(x_k)$.
\end{proof}
\end{lemma}

\begin{corollary}\label{11_}
	
There is no AS-Ore extension with Hilbert series
\[h(t)=\frac {1} {(1-t)^{k}\prod \limits _{j=1}^{n-k} (1-t^{i_{j}})},\;  i_{j}>2\; \text{for all } j,\; k<n.
\] 

\begin{proof}
Suppose this is possible.  Choose $k$ such that $deg(x_k)$ is minimal amongst variables with degree greater than 1.  By \Cref{ijk}, $deg(x_k)=deg(x_i)+deg(x_j)$ for some $i$ and $j$.  If $deg(x_i)=deg(x_j)=1$ then this equation says $deg(x_k)=2$, but an Ore extension with the given the Hilbert series cannot have any variables of degree two.  Otherwise, we can assume that $deg(x_i)>1$.  Since $x_k$ was chosen to have smallest degree greater than 1, the equation now becomes $deg(x_k)=deg(x_i)+deg(x_j)\geq deg(x_k)+1$, which is impossible.  Thus, no Ore extension generated in degree one can have this Hilbert series.
\end{proof}
\end{corollary}

\iffalse
\begin{proposition}
There is no iterated Ore extension generated in degree one with Hilbert series 
\[h(t)=\frac {1} {(1-t)^{2}\prod \limits _{j=1}^{n-2} (1-t^{i_{j}})},\;  i_{j}>2\; \forall j
\]

\begin{proof}
Assume the contrary.  Then there are two degree one variables, call them $x_{n}$ and $x_{n-1}$, and other variables of degree at least 3. Then $x_{n}x_{n-1}=\sigma (x_{n-1})x_{n}+\delta (x_{n-1})$ is a degree two homogeneous polynomial with leading term $x_{n}x_{n-1}$ so must have the form $ x_{n}x_{n-1}=a_{1}x_{n-1}x_{n}+a_{2}x_{n-1}^{2}$ and so the subring generated by $x_{n}$ and $x_{n-1}$, call it $K\langle x_{n-1},x_{n}\rangle$, has basis $\{ x_{n-1}^{i}x_{n}^{j}\; | \; i,j\in \mathbb Z_{\ge 0} \}$.
 
If $x_{i} \in K\langle x_{n-1},x_{n}\rangle$,\; $i<n-1$, then 
\[x_{i}=\sum \limits_{k=0}^{d_{i}} a_{k}x_{n-1}^{k}x_{n}^{d_{i}-k}, \; d_{i}=deg(x_{i})
\]
violates the condition for being an Ore extension since either $x_{i}=0$ or this relation has leading term $a_{j}x_{n-1}^{j}x_{n}^{d_{i}-k}$ where $j=min\{k\; |\; a_{k}\neq 0\}$ while all leading terms of the Ore extension should have the form $x_{j}x_{i}$, $j>i$.
\end{proof}

\end{proposition}

\fi

\begin{lemma}\label{1122}
There is no AS-Ore extension with Hilbert series 
\[
h(t)=\frac {1} {(1-t)^{2}(1-t^{2})^{2}\prod \limits _{j=1}^{n-4} (1-t^{i_{j}})},\;  i_{j}\geq 2\; \text{for all } j.
\]

\begin{proof}
To find a contradiction, assume that such an extension, $A$, exists.  By \Cref{Presentation 1}, there exists an ordering on the variables such that $A\cong \displaystyle \frac{K\langle x_1,\cdots, x_n\rangle}{\langle \{r_{ji} \} \rangle}$. Let $x_{n}$ and $x_{n-1}$ be degree one variables and let $x_{n-2}$ and $x_{n-3}$ be distinct degree two variables in $K\langle x_{n-1},x_{n}\rangle$ and so of the form 
\begin{align*}
& x_{n-2}=a_{1}x_{n-1}^2+a_{2}x_{n-1}x_{n}+a_{3}x_{n}x_{n-1}+a_{4}x_{n}^{2},\\
& x_{n-3}=b_{1}x_{n-1}^2+b_{2}x_{n-1}x_{n}+b_{3}x_{n}x_{n-1}+b_{4}x_{n}^{2}.
\end{align*}
Without loss of generality, we can assume that $x_n>x_{n-1}$ in the ordering.  Also, we must have that $x_{n-2}<x_n$ and $x_{n-3}<x_n$ since an Ore extension cannot have $x_{n-2}$ or $x_{n-3}$ as the leading term of a relation.  Since an Ore extension has no leading term of the form $x_{n}^{2}$ and a unique term of the form $x_{n}x_{n-1}$, we get that $a_{4}=b_{4}=0$ and one of $a_{3}$ and $ b_{3}$ is 0. Without loss of generality, assume that $b_{3}=0$. Then the relation $x_{n-3}=b_{1}x_{n-1}^2+b_{2}x_{n-1}x_{n}$ has a leading term inconsistent with an Ore extension.
\end{proof}
\end{lemma}

\iffalse
In particular, there is no Ore extension generated in degree one with Hilbert series
\[
h(t)=\frac {1} {(1-t)^{2}(1-t^{2})^{2}(1-t^{3})}
\]
\fi

\begin{theorem} \label{DegreeType}
For an AS-Ore extension with 5 variables, one of the following options represents the possible degree type of the extension.
\begin{enumerate}
\item $(1,1,2,3,5),$
\item $(1,1,2,3,4),$
\item $(1,1,2,3,3),$
\item $(1,1,1,2,3),$
\item $(1,1,1,2,2),$
\item $(1,1,1,1,2),$
\item $(1,1,1,1,1).$

\end{enumerate}

\begin{proof}
Clearly an Ore extension with no variables of degree one cannot be generated in degree one.

Similarly there is no Ore extension generated in degree one with just 1 degree one variable.  For if there were such an Ore extension and $x_k$ were of minimal degree amongst the remaining 4 variables, \Cref{ijk} says that $deg(x_k)=deg(x_i)+deg(x_j)\geq1+deg(x_k)$ 
and such an inequality is impossible.

If the Ore extension has exactly 2 degree one variables then \Cref{11_} implies that there is at least one variable of degree two.  If $x_k$ is of minimal degree amongst the remaining 2 variables then \Cref{ijk} tells us that $deg(x_k)\in \{2,3\}$ since these are the only possible combinations of $deg(x_i)+deg(x_j)$.  By \Cref{1122}, there can be at most one variable of degree two, so $deg(x_k)=3$.  Again by \Cref{ijk}, the final and largest degree variable, $x_l$, satisfies $deg(x_l)=deg(x_i)+deg(x_j)$ for some $i$ and $j$ so $deg(x_l)\in \{3,4,5\}$.  Thus, the list of possible degree types for an Ore extension with exactly 2 degree one variables is
\begin{enumerate}
\item $(1,1,2,3,5),$
\item $(1,1,2,3,4),$
\item $(1,1,2,3,3).$
\end{enumerate}

If the Ore extension has exactly 3 degree one variables then it must have at least one degree two variable and the remaining variable, by \Cref{ijk}, must be of degree two or three.  Thus, the list of possible degree types in this case is
\begin{enumerate}
\setcounter{enumi}{3}
\item $(1,1,1,2,3),$
\item $(1,1,1,2,2).$
\end{enumerate}

If the Ore extension has exactly 4 degree one variables, then \Cref{ijk} tells us that the remaining variable must be degree two and the possible degree type is
\begin{enumerate}
\setcounter{enumi}{5}
\item $(1,1,1,1,2).$
\end{enumerate}

Finally, it is possible for the Ore extension to have 5 degree one variables and degree type
\begin{enumerate}
\setcounter{enumi}{6}
\item $(1,1,1,1,1).$
\end{enumerate}
\vspace{-.5cm}
\end{proof}
\end{theorem}

%\vspace{-.5 cm}
While this result technically only restricts the possible degree types of AS-Ore extensions, it is also true that, for each of the 7 possible options listed, there exists an AS-Ore extension with the given type.  A commutative ring in five variables is an example of an algebra with type (1,1,1,1,1).  For options 2-6, there are enveloping algebras with variables of appropriate degrees (see \cite[Section 3]{FV},  \Cref{EnvelopingExamples}, \Cref{Env(11122)}, \Cref{Env(11123)}). Finally, we construct an AS-Ore extension with degree type (1,1,2,3,5) in the next section (\Cref{Ore11235}).

%%%%%%%%%%%%%%%%%%%%%%%%%%%%%%%%%%%%%%%%%%%%%%%%%%%%%%%%%%%%%%%%%%%%%%%%%%%%

\section{An AS-Ore extension with degree type $(1,1,2,3,5)$}

For AS-regular algebras of dimension at most 4, it is known that every Hilbert series has a unique relation type and every relation type can be realized by the enveloping algebra of a graded Lie algebra. In their paper, Floystad and Vatne asked whether this held in dimension 5 and constructed, as a counter example, an AS-regular algebra with 2 degree one generators and Hilbert series $\displaystyle h(t)= \frac {1} {(1-t)^{2}(1-t^{2})(1-t^{3})(1-t^{5})}$.  By looking at the shifts in the free resolution, they prove that there is no enveloping algebra of a graded Lie algebra with this Hilbert series \cite[Proposition 3.4 and Theorem 4.2]{FV}.  Based on the presentation of an enveloping algebra given above, we provide an alternate proof of this result.

\begin{proposition}
There is no enveloping algebra of a graded Lie algebra which is generated in degree one with  $\displaystyle  h(t)= \frac {1} {(1-t)^{2}(1-t^{2})(1-t^{3})(1-t^{5})}$.

\begin{proof}
Assume to the contrary that there is such an algebra.  Then there are variables $(x_5,x_4,x_3,x_2,x_1)$ with respective degrees $(1,1,2,3,5)$ . Consider the possible terms in the relations.  We will list only those required to show the contradiction.
\begin{align*}
r_{54}: x_5x_4&=x_4x_5+a_1x_3\\
r_{52}: x_5x_2&=x_2x_5\\
r_{42}: x_4x_2&=x_2x_4\\
r_{32}: x_3x_2&=x_2x_3+b_1x_1.
\end{align*}
Here $a_1$ and $ b_1$ must be nonzero for this algebra to be generated in degree one (by the proof of \Cref{ijk}) since these are the only relations of degree two and five respectively.  Additionally, the middle two relations cannot contain any additional terms since they are degree four and this algebra contains no variable of degree exactly 4.
Now consider:
\begin{align*}
x_5(x_4x_2)&=x_5x_2x_4\\
&=x_2x_5x_4\\
&=x_2x_4x_5+a_1x_2x_3, \text{ while} \\
(x_5x_4)x_2&=x_4x_5x_2+a_1x_3x_2\\
&=x_4x_2x_5+a_1x_2x_3+a_1b_1x_1\\
&=x_2x_4x_5+a_1x_2x_3+a_1b_1x_1.
\end{align*}
In order for this overlap to resolve, $a_1b_1=0$, which is impossible if this algebra is generated in degree one.
\end{proof}
\end{proposition}

 It is natural to ask whether every relation type can be realized by a generalization of an enveloping algebra, in particular an AS-Ore extension.  In the literature for algebras of dimension 5 there are currently only two known relation types that cannot be realized by an enveloping algebra.  One has Hilbert series $\displaystyle h(t)=\frac{1}{(1-t)^2(1-t^2)(1-t^3)^2}$, relation type $(4,4,4,5,5)$, and can be realized by an AS-Ore extension \cite[Section 5.2]{WW}.

In support of the hypothesis that all relation types can be realized by an AS-Ore extension, we present an example of an AS-Ore extension with the other known relation type that cannot be realized by an enveloping algebra.  This is equivalent to finding an example of an AS-Ore extension with the appropriate Hilbert series since Floystad and Vatne have already classified the possible relation types of algebras with two generators \cite[Theorem 5.6]{FV}, and the relation type of an algebra with this Hilbert series is unique.

The following example was found by writing the general relations provided by \Cref{Presentation 1} and using the mathematical software program Mathematica to solve the large system of equations that result from setting overlap ambiguities equal to 0.  This proved to be an overwhelming project for the computer and some coefficients were ultimately assumed to be 0 to make the computations possible, as our goal was to prove the existence of such an algebra rather than to completely classify the possible families of relations.

\begin{theorem}\label{Ore11235}
The following relations define an AS-Ore extension (an iterated Ore extension which is graded, generated in degree one, and has each $\sigma_i$ an injection)  which has\\ $\displaystyle  h(t)= \frac {1} {(1-t)^{2}(1-t^{2})(1-t^{3})(1-t^{5})}$ and relation type (3,4,7):
\begin{align*}
r_{21}:x_{2}x_{1}&=-x_{1}x_{2}\\
r_{32}:x_{3}x_{2}&=x_{1}+b x_{2}x_{3}\\
r_{31}:x_{3}x_{1}&=-x_{1}x_{3}\\
r_{43}:x_{4}x_{3}&=x_{2}+bx_{3}x_{4}\\
r_{42}:x_4x_2&=b^2x_2x_4\\
r_{41}:x_{4}x_{1}&=x_{1}x_{4}\\
r_{54}:x_{5}x_{4}&=x_{3}+x_{4}x_{5}\\
r_{53}:x_{5}x_{3}&=-x_{3}x_{5}\\
r_{52}:x_5x_2&=-x_2x_5-b^2x_3x_3\\
r_{51}:x_5x_1&=x_1x_5+cx_3x_3x_3,
\end{align*}
where  $\displaystyle b=e^{\frac{4 \pi i}{3}}$ \text{and}  $c=\frac{2b^2}{1-b+b^2}$.

\begin{proof}

Let the degrees of $(x_5,x_4,x_3,x_2,x_1)$ be $(1,1,2,3,5)$ with lexicographic order $x_5>\cdots>x_1$ so that the leading terms are as presented. 

We note that, for the given degrees of these variables, each of the above relations is homogeneous.  To check that this is an Ore extension, we then check that all reduction ambiguities resolve.  All ambiguities have the form $x_kx_jx_i$ where $k>j>i$ and there are a total of 10 such ambiguities for this set of relations.  A computation shows that, for the given choice of $b$ and $c$, all overlaps resolve.  We carry out this computation in Mathematica and the code for these and future computations can be found on the author's website \cite[Section 1]{SE}.  Thus, $\{x_1^{e_1} \cdots x_n^{e_n}\}$ is a basis for this algebra and this is Ore by \Cref{Presentation2}. 

It remains to check that this is generated in degree one and that $\sigma_j$ is injective for all j.  To see that this algebra is generated in degree one, note that $r_{32}$, $r_{43}$, and $r_{54}$ can be solved for $x_1$, $x_2$, and $x_3$ respectively and so everything may be expressed in terms of the degree one generators $x_4$ and $x_5$.
\iffalse  So we have
\begin{align*}
x_3&=x_5x_4-x_4x_5,\\
x_2&=x_4x_3-b^2x_3x_4\\
&=-bx_5x_5x_4+(1+b)x_4x_5x_4-x_4x_4x_5,\\
x_1&=x_3x_2-bx_2x_3\\
&=-b x_5x_4x_5x_4x_4+(1+b+b^2)x_5x_4x_4x_5x_4+(-1-b^2)x_5x_4x_4x_4x_5+bx_4x_5x_5x_4x_4+\\
& \; \; \; (-1-2b-b^2)x_4x_5x_4x_5x_4+(1+b+b^2)x_4x_5x_4x_4x_5+bx_4x_4x_5x_5x_4-bx_4x_4x_5x_4x_5.
\end{align*}
Since $x_4$ and $x_5$ are both degree one and the other generators can be written in terms of these variables, the algebra is indeed generated in degree one.
\fi
Note that for all $1\leq i<j\leq 5$, $\sigma_j(x_i)=a_{ji}x_i$ where $a_{ji}$ is a root of unity.  Thus, $\sigma_j^{n_j}$ is the identity map for some $n_j$ and so each homomorphism is injective.  Thus, this algebra is an AS-Ore extension.

That the algebra has the desired Hilbert series is now immediate from the fact that is has the same basis, and therefore the same Hilbert series, as the weighted commutative polynomial ring with the same variables and degrees.  We again note that the relation type of an algebra with this Hilbert series is known to be (3,4,7) \cite[Theorem 5.6]{FV}, although the computations proving it in this case are also included in the online code.
\end{proof}

\end{theorem}

\section{A classification of relation types for AS-Ore extensions with 4 degree one generators}

We now begin the process of attempting to classify all possible resolution types of dimension 5 AS-Ore extensions generated in degree one, beginning with the case where the algebra has 4 degree one generators.  It will be convenient to alternate between thinking of the algebra as an Ore extension with the presentation given by \Cref{Presentation 1}, $\displaystyle A\cong \frac{K\langle x_1,\cdots,x_n\rangle}{\langle\{r_{ji}\}\rangle}$, and as an algebra presented in terms of its degree one generators, $\displaystyle A\cong\A\text:=\frac{K\langle x_1,\cdots,x_b\rangle}{I}$.  We will fix the notation that $A$ refers to the algebra viewed as an AS-Ore extension presented by 5 generators and that $\A$ is an algebra isomorphic to $A$, viewed as generated in degree one.  We get it from $A$ by changing the ordering on the variables so that $x_i>x_j$ whenever $\deg(x_i)>1$ and $\deg(x_j)=1$, making a choice that allows us to solve the Ore relations for all variables that are not degree one, and writing the remaining relations in terms of the degree one generators.  Changing the ordering of the variables may change the \basis $ $ of the algebra, but will of course not change the Hilbert series.  By the construction of $\A$, we also note that its minimal generating set cannot contain more elements of a particular degree than what the minimal generating set of $A$ has.

 By \Cref{Resolution}, the free resolution of any dimension 5 regular algebra generated in degree one is
\[0\rightarrow A(-l)\rightarrow A(-l+1)^b\rightarrow \bigoplus \limits_{i=1}^{n} A(-l+a_i)\rightarrow \bigoplus \limits_{i=1}^{n}  A(-a_i)\rightarrow A(-1)^b\rightarrow A\rightarrow K\rightarrow 0\] where $b$ represents the number of degree one generators, $l$ the total shift of the resolution, $n$ the number of relations in the minimal generating set of the ideal $I$, and $a_i$ the homogeneous degree of the $i$th relation of a fixed minimal generating set of $I$.  By the symmetry of this free resolution, the resolution type is uniquely determined by the $a_i$ together with the Hilbert series of the algebra since the series will determine the value of $l$.  Thus, it suffices to classify the possible relation types $(a_1,\cdots,a_n)$, $a_i\leq\cdots\leq a_n$, of dimension 5 AS-Ore algebras.  
%We note that this resolution is only valid for algebras generated in degree one, so we will alternate between thinking about an Ore extension, $A$, which is not generated in degree one, and the algebra $\A$, which is isomorphic to $A$ and generated in degree one.
%will have to consider an ordering on the variables that allows us to view our Ore extension as an algebra generated in degree one, and examine the minimal generating set of the ideal defined by relations under this new ordering.

Recall by \Cref{DegreeType} that an AS-Ore extension, $A$, with 5 variables and 4 degree one generators will have degree type (1,1,1,1,2).  Since the relations for an Ore extension come from the $\{ r_{ji} \}$ as described in \Cref{Presentation 1} (and there are no additional relations since overlaps resolve by the same theorem), $A$ will have 6 degree two relations in the \basis $ $, 4 relations of degree three, and no relations of degree four or larger.  

Let $\widetilde{A}$ denote the same algebra, $A$, viewed as an algebra generated in degree one via the process explained above where 1 degree two relation of ${A}$ will be used to express the degree two variable in terms of the generators and the remaining 5 will be part of the minimal generating set of $\A$.  Since $A$ has 4 degree three relations, $\A$ will have at most 4 degree three relations.  It is possible that $\widetilde{A}$ will have fewer than 4 $K$-independent degree three relations or for these relations to be consequences of overlaps that fail to resolve rather than part of the minimal generating set.  Since $A$ has no relations of degree more than three, any relations of degree more than three in $\A$ must be consequences of overlaps that fail to resolve and so not part of the minimal generating set of $\widetilde{A}$.  Thus, the only candidates for the relation type of an AS-Ore extension with degree type (1,1,1,1,2) are:
\begin{enumerate}
	\item (2,2,2,2,2),
	\item (2,2,2,2,2,3),
	\item (2,2,2,2,2,3,3),
	\item(2,2,2,2,2,3,3,3),
	\item(2,2,2,2,2,3,3,3,3).
\end{enumerate}

In the next theorems, we will classify all possible relation types of an AS-Ore extension with the given degree type.  We prove that types (4) and (5) above are impossible, that types (1) and (3) can be realized by enveloping algebras, and that type (2) can be realized by an AS-Ore extension but not by an enveloping algebra.  This differs slightly from the comment in \cite[Section 3]{FV} where examples of enveloping algebras of types (1) and (3) are explicitly presented but the reader is encouraged to also check that there is an example of an enveloping algebras of type (2). \\

In order to prove that certain relation types are impossible, we need to know more about the specific leading terms and the overlaps that come from the degree two relations.  We use a simplified version of Hilbert driven Gr{\"o}bner basis computation, a technique used by Rogalski and Zhang to study $\mathds{Z} ^2$-graded dimension 4 algebras with 3 generators \cite{RZ} and later used by Zhou and Lu to classify possible families of relations of  $\mathds{Z} ^2$-graded dimension 5 algebras with 2 generators \cite{ZL}.

Let $A$ be an AS-Ore extension with degree type (1,1,1,1,2).  Then $A$ has Hilbert series $\displaystyle h_A(t)=\frac{1}{(1-t)^4(1-t^2)}=1+4t+11t^2+24t^3+O(t^4)$.  The idea of Hilbert driven basis computation is to construct $\A$ by viewing it as a free algebra on its degree one  generators modulo an ideal $I$, and to identify the generators of $I$ by comparing the Hilbert series of the constructed algebra against the known Hilbert series.  In order to do this one dimension at a time, we will use a \emph{monomial algebra} which we get by replacing each relation in the \basis $ $ with just the leading term of the relation.  The details of this construction can be found in \cite[Section 2]{ZL}, along with the proof that the monomial algebra will have the same Hilbert series as the original algebra \cite[Lemma 2.1]{ZL}.  Let $\A_0$ denote the free algebra on four generators.  Then $h_{\A_0}(t)=1+4t+16t^2+64t^3+O(t^4)$ and $h_{\A_0}(t)-h_{\A} (t)=5t^2+O(t^3)$.  Thus, $I$ must contain 5 degree two relations.  Although we already knew this, the method can be used to find the number of relations in the basis of higher degrees, although the analysis does depend on which leading terms we choose for our relations. 

Let $x_1$ be the degree two variable and without loss of generality, list the degree one variables so that $x_2<x_3<x_4<x_5$.  The degree two LT's (leading terms) in $\A$ must come from the $\{r_{ji}$\} relations of the Ore extension, $A$, and so must belong to the list  $\{x_3x_2,x_4x_3, x_4x_2, x_5x_4, x_5x_3, x_5x_4\}$. If the set of degree two LTs in $I$ is \linebreak  $\{x_3x_2,x_4x_3, x_4x_2, x_5x_4, x_5x_3\}$, then let us denote the monomial algebra which has these leading terms as its relations by $\A_2$ since the Hilbert series agrees with that of $\A$ up to dimension 2.  Then $h_{\A_2}-h_{\A}=4t^3+O(t^4)$ so there must be 4 degree three relations in the Gr{\"o}bner basis.  The calculations for the difference of these Hilbert series were done in Mathematica, and the code is available online \cite[Section 2]{SE}.

If instead we start with LTs in $I$ $\{x_3x_2,x_4x_3,  x_5x_4, x_5x_3,x_5x_2\}$ or \linebreak $\{x_3x_2,x_4x_3, x_4x_2, x_5x_4, x_5x_2\}$, then $h_{\A_2}-h_{\A}=3t^3+O(t^4)$ and there are 3 degree three relations in the basis.  If we start with LTs $\{x_4x_3, x_4x_2, x_5x_4, x_5x_3,x_5x_2\}$, $\{x_3x_2, x_4x_2, x_5x_4, x_5x_3,x_5x_2\}$, or $\{x_3x_2,x_4x_3, x_4x_2, x_5x_3,x_5x_2\}$ then there will be 2 degree three relations by a similar analysis.

So we have found that the Gr{\"o}bner basis of an AS-Ore extension with the given degree type has 5 degree two relations and 2-4 degree three relations.  We could continue the process to see what the possible LTs of degree three are and if the basis has additional relations of higher degree, but it is not useful for our analysis.  We remain more interested in the number and degrees of the minimal generators since these completely classify the possible relation types of the algebra, and we know the algebra has no minimal relations of degree greater than three.  We still need to investigate, in each case, which of the degree three relations are part of the minimal generating set and which are simply consequences of overlaps that fail to resolve.

\begin{theorem} \label{No(3333)}
	There is no AS-Ore extension with degree type $(1,1,1,1,2)$ and minimal relation type $(2,2,2,2,2,3,3,3,3)$.
	
	\begin{proof}
		Assume to the contrary that there is such an AS-Ore extension \linebreak $A=K[x_{i_1}]\cdots [x_{i_5},\sigma_{i_5},\delta_{i_5}]$.  Label the degree two variable $x_1$ and label the degree one variables so that $x_2<x_3<x_4<x_5$ in the ordering. (We make no assumption about when $x_1$ is adjoined.) The list of reduced degree 2 monomials is 
		\[\{x_ 1, x_ 2 x_ 2, x_ 2 x_ 3, x_ 3 x_ 3, x_ 2 x_ 4, x_ 3 x_ 4, x_ 4 x_ 4, x_ 2 x_ 5, x_ 3 x_ 5, x_ 4 x_ 5, x_ 5 x_ 5\}.\]
		From \Cref{Presentation 1}, $x_j$ must occur only to the first power in the relation with leading term $x_jx_i$ and $x_k$ should not appear for any $x_k>x_j$.  Based on these observations, we will write the most general possible degree two relations:
		\begin{align*}
		r_{32}:x_ 3 x_ 2 &= b_ 1 x_ 1 + b_ 2 x_ 2 x_ 2 + b_ 3 x_ 2 x_ 3\\
		r_{42}:x_ 4 x_ 2 &=  e_ 1 x_ 1 + e_ 2 x_ 2 x_ 2 + e_ 3 x_ 2 x_ 3 + e_ 4 x_ 2 x_ 4 +  e_ 6 x_ 3 x_ 3 + e_ 7 x_ 3 x_ 4\\
		r_{43}:x_ 4 x_ 3 &=  d_ 1 x_ 1 + d_ 2 x_ 2 x_ 2 + d_ 3 x_ 2 x_ 3 + d_ 4 x_ 2 x_ 4 +   d_ 6 x_ 3 x_ 3 + d_ 7 x_ 3 x_ 4\\
		r_{52}:x_ 5 x_ 2 &= i_ 1 x_ 1 + i_ 2 x_ 2 x_ 2 + i_ 3 x_ 2 x_ 3 + i_ 4 x_ 2 x_ 4 +   i_ 5 x_ 2 x_ 5 + i_ 6 x_ 3 x_ 3 + i_ 7 x_ 3 x_ 4 + i_ 8 x_ 3 x_ 5 \\
		& +   i_ 9 x_ 4 x_ 4 + i_ {10} x_ 4 x_ 5\\
		r_{53}:x_ 5 x_ 3 &=  h_ 1 x_ 1 + h_ 2 x_ 2 x_ 2 + h_ 3 x_ 2 x_ 3 + h_ 4 x_ 2 x_ 4  +   h_ 5 x_ 2 x_ 5 + h_ 6 x_ 3 x_ 3 + h_ 7 x_ 3 x_ 4 \\
		& + h_ 8 x_ 3 x_ 5 +   h_ 9 x_ 4 x_ 4 + h_ {10} x_ 4 x_ 5\\
		r_{54}:x_ 5 x_ 4 &=  g_ 1 x_ 1 + g_ 2 x_ 2 x_ 2 + g_ 3 x_ 2 x_ 3 + g_ 4 x_ 2 x_ 4 +   g_ 5 x_ 2 x_ 5 + g_ 6 x_ 3 x_ 3 + g_ 7 x_ 3 x_ 4 \\
		& 
		+ g_ 8 x_ 3 x_ 5 +   g_ 9 x_ 4 x_ 4 + g_ {10} x_ 4 x_ 5.
		\end{align*}
		(It is worth noting that if $x_1$ is adjoined late, some of these coefficients must be zero, although this fact will not be needed to complete the contradiction.  For example, if the Ore extension is $K[x_2][x_3,\sigma_3,\delta_3][x_1,\sigma_1,\delta_1][x_4,\sigma_4,\delta_4][x_5,\sigma_5,\delta_5]$, then $b_1$ must be zero since $x_3x_2=\sigma_3(x_2)x_3+\delta_3(x_2)$ where $\delta_3$ is a derivation of $K[x_2]$ and thus cannot map $x_2$ to a term containing $x_1$.)
		
		Without loss of generality, we can solve for $x_1$ using the relation with smallest leading term that has a nonzero coefficient of $x_1$.  For example, if $b_1\neq0$ then we can solve $r_{32}$ to find that $x_1=\displaystyle \frac{1}{b_1}(x_3x_2-b_2x_2x_2-b_3x_2x_3)$.  Consider now the \basis $ $ of the algebra $\A$, viewed as an algebra generated in degree one.
		
		The preceding analysis shows that the list of degree two LTs of $\A$ must be \linebreak $\{x_3x_2,x_4x_3, x_4x_2, x_5x_4, x_5x_3\}$ since this is the only set of LTs that has 4 degree three relations in the \basis.  In particular, $x_5x_2$ is not a leading term.  But if any of $b_1$, $e_1$, or $d_1$ are non-zero, then $x_1$ can be expressed in terms of monomials smaller than $x_5x_2$ and the LT of $r_{52}$ will be $x_5x_2$, a contradiction.  We may therefore assume that $b_1=e_1 =d_1=0$ and we may assume that $i_1$ is not zero since again this would otherwise make the leading term of $r_{52}$ equal to $x_5x_2$. Thus $\A$ is obtained by using $r_{52}$ to write $\displaystyle x_1=\frac{1}{i_1}(x_5x_2-i_{10}x_4x_5-\cdots -i_2x_2x_2)$ and substituting this expression for $x_1$ in the other relations. 
		
		We are interested in the coefficient of $x_5x_2x_3$ in the reduction of the overlap $x_5x_3x_2$ in $\widetilde{A}$ since this will provide the contradiction. We compute:
		\begin{align*}
		x_5(x_3x_2)&=x_5(0x_1+b_3x_2x_3+[\text{smaller terms}])\\
		&=b_3x_5x_2x_3+[\text{smaller terms}], \text{ while}\\
		(x_5x_3)x_2&=(h_1x_1+h_{10}x_4x_5+[\text{smaller terms}])x_2\\
		&=(\frac{h_1}{i_1}(x_5x_2-i_{10}x_4x_5-[\text{smaller terms}])+h_{10}x_4x_5+[\text{smaller terms}])x_2\\
		&=\frac{h_1}{i_1}x_5x_2x_2+[\text{smaller terms}].
		\end{align*}
	Note that $x_5x_2x_3$ is a reduced word with respect to the degree two relations of $\A$.  Thus $x_5(x_3x_2)-(x_5x_3)x_2=b_3x_5x_2x_3+[\text{smaller terms}]$ and $b_3=0$ if this overlap resolves.  However, if $b_3=0$ then $r_{32}$ becomes
		\[x_3x_2=\sigma_3(x_2)x_3+\delta_3(x_2)=0x_2x_3+0x_1+b_2x_2x_2.\]
		This would suggest that $\sigma_3(x_2)=0$ and $\sigma_3$ is not injective, which contradicts the claim that the initial algebra was an AS-Ore extension.  Thus, $b_3$ is not zero, this overlap does not resolve, at least one of the degree three relations in the 
		\basis $ $ of this algebra is a consequence of an overlap, and so there are not 4 degree three relations in the minimal generating set.		
	\end{proof}
\end{theorem}

With a little more work, the same technique can be used to show that there cannot be 3 degree three relations in the minimal generating set.

\begin{theorem} \label{No(333)}
	There is no AS-Ore extension generated in degree one with degree type $(1,1,1,1,2)$ and relation type $(2,2,2,2,2,3,3,3)$.
	\begin{proof}
		Assume to the contrary.  As in the previous theorem, label the degree two variable $x_1$ and label the degree one variables so that $x_2<x_3<x_4<x_5$ in the ordering.  Then the general form of the possible degree two relations is the same as in the previous theorem.  We will consider the different cases where the set of degree three leading terms leads to Gr{\"o}bner bases with at least 3 degree three relations.\\
		
		\textbf{Case 1:} The set of degree two LTs in $\widetilde{A}$ is $\{x_3x_2,x_4x_3, x_4x_2, x_5x_4, x_5x_3\}$.\\
		By the same argument as that in \Cref{No(3333)}, since $x_5x_2$ is not a leading term, the relations with smaller leading terms must have coefficient in front of $x_1$ equal to zero and $r_{52}$ must have a nonzero coefficient in front of the $x_1$.  Thus, $b_1=e_1=d_1=0$, $i_1\neq0$, there are 4 degree three relations in the \basis, and we have already seen in the proof of \Cref{No(333)} that $x_5x_3x_2$ is an ambiguity that fails to resolve.  Taking $x_5(x_3x_2)-(x_5x_3)x_2$ in $\widetilde{A}$ gives a new degree three relation with leading term $x_5x_2x_3$.  When evaluating whether other degree three overlaps resolve, we should reduce them modulo this new relation, but this will not be necessary in the calculations that follow since $x_5x_2x_3$ is too small to affect the analysis of whether the overlaps resolve.  Now consider the coefficient of $x_5x_2x_4$ in the following reduction of overlaps in $\widetilde{A}$:
		\begin{align*}
		x_5(x_4x_2)&=x_5(0x_1+e_4x_2x_4+e_6x_3x_3+e_7x_3x_4+[\text{smaller terms}])\\
		&=e_4x_5x_2x_4+e_6x_5x_3x_3+e_7x_5x_3x_4+[\text{small}]\\
		&=e_4x_5x_2x_4+e_6(h_1x_1+[\text{small}])x_3+e_7(h_1x_1+[\text{small}])x_4 +[\text{small}]\\
		&=e_4x_5x_2x_4+e_6(\frac{h_1}{i_1}x_5x_2-[\text{small}])x_3+e_7(\frac{h_1}{i_1}x_5x_2-[\text{small}])x_4+[\text{small}]\\
		&=(e_4+\frac{e_7h_1}{i_1})x_5x_2x_4+[\text{small}], \text{while}\\
		(x_5x_4)x_2&=(g_1x_1+[\text{smaller terms}])x_2\\
		&=\frac{g_1}{i_1}x_5x_2x_2+[\text{small}]
		\end{align*}
		So $\displaystyle x_5(x_4x_2)-(x_5x_4)x_2=(e_4+\frac{e_7h_1}{i_1})x_5x_2x_4+[\text{smaller terms}]$.
		
	Similarly, $\displaystyle x_5(x_4x_3)-(x_5x_4)x_3=(d_4+\frac{d_7h_1}{i_1})x_5x_2x_4+[\text{smaller terms}]$.
		From the relations 
		\begin{align*}
		r_{42}:x_ 4 x_ 2 &=  e_ 1 x_ 1 + e_ 2 x_ 2 x_ 2 + e_ 3 x_ 2 x_ 3 + e_ 4 x_ 2 x_ 4 +  e_ 6 x_ 3 x_ 3 + e_ 7 x_ 3 x_ 4 \text{ and}\\
		r_{43}:x_ 4 x_ 3 &=  d_ 1 x_ 1 + d_ 2 x_ 2 x_ 2 + d_ 3 x_ 2 x_ 3 + d_ 4 x_ 2 x_ 4 +   d_ 6 x_ 3 x_ 3 + d_ 7 x_ 3 x_ 4,
		\end{align*}
		we see that 
		\begin{align*}
		\sigma_4(x_2)&=e_4x_2+e_7x_3 \text{ and}\\
		\sigma_4(x_3)&=d_4x_2+d_7x_3.
		\end{align*}
		Thus, $\sigma_4$ is injective if and only if 
		$\displaystyle \det \begin{bmatrix}
		e_4 &e_7\\d_4 &d_7
		\end{bmatrix}\neq 0.$
		If we assume that both overlaps resolve, $\displaystyle e_4=-\frac{e_7 h_1}{i_1}$ and $\displaystyle d_4=-\frac{d_7 h_1}{i_1}$, so 
		$  \begin{vmatrix}
		e_4 &e_7\\d_4 &d_7
		\end{vmatrix}\ =\begin{vmatrix}
		\displaystyle \frac{-e_7 h_1}{i_1}& e_7\\\displaystyle \frac{-d_7 h_1}{i_1} &d_7
		\end{vmatrix}=0$ and $\sigma_4$ is not injective, a contradiction.  Thus, at least one of the two overlaps, $x_5x_4x_3$ or $x_5x_4x_2$, fails to resolve.  In total, there are at least two overlaps that do not resolve so there are at most 2 degree three relations in the minimal generating set.\\
		
		\iffalse
		and 
		\begin{align*}
		\sigma_4(d_4x_2-e_4x_3)&=d_4e_4x_2+d_4e_7x_3-e_4d_4x_2-e_4d_7x_3\\
		&=(d_4e_4-e_4d_4)x_2+(d_4e_7-e_4d_7)x_3\\
		&=(d_4e_4-e_4d_4)x_2+(-\frac{d_7h_1e_7}{i_1}+\frac{e_7h_1d_7}{i_1})x_3\\
		&=0.
		\end{align*}
		This violates the fact that $\sigma_4$ is injective unless $e_4$ and $d_4$ are both zero.  But in this case, $\sigma_4(d_7x_2-e_7x_3)=0$ and $\sigma_4$ is not injective unless $d_7=0$ and $e_7=0$.  If $e_7=0$ then $\sigma_4(x_2)=0$.
		\fi

		\textbf{Case 2:} The set of degree two LTs in $\widetilde{A}$ is $\{x_3x_2,x_4x_3, x_4x_2, x_5x_4, x_5x_2\}$.\\
		In this case, $x_5x_3$ is not a leading term and similar reasoning as that used in  \Cref{No(3333)} allows us to conclude that $b_1=e_1=d_1=i_1=0$, $h_1\neq0$, and there are 3 degree three relations in the Gr{\"o}bner basis. We then compute $x_5(x_4x_2)-(x_5x_4)x_2=e_7 x_5x_3x_4+ [\text{smaller terms}]$ and $x_5(x_4x_3)-(x_5x_4)x_3=d_7 x_5x_3x_4+ [\text{smaller terms}]$.  These computations can be done by hand by looking at the largest terms in the reduction, just as in the previous case, but we omit the details.  The code used for all calculations in this proof is on the author's website \cite[Section 3]{SE}.  If these overlaps both resolve, $e_7=0$ and $d_7=0,$ $\begin{vmatrix} e_4 &e_7\\d_4 &d_7 \end{vmatrix}=0,$ and $\sigma_4$ is not injective, which is a contradiction.  Thus, one of these overlaps must fail to resolve and the minimal generating set has at most 2 relations of degree three.\\
		
		\textbf{Case 3:} The set of degree two LTs in $\widetilde{A}$ is $\{x_3x_2,x_4x_3, x_5x_4, x_5x_3,x_5x_2\}$.\\
		In this case, $x_4x_2$ is not a leading term so $b_1=0$, $e_1\neq0$, and there are 3 degree three relations in the \basis. Then $x_4(x_3x_2)-(x_4x_3)x_2=b_3 x_4x_2x_3+\text{[smaller terms]}$.  This computation can be done by hand, but we omit the details here.  If $\sigma_3$ is injective, $b_3$ cannot be zero.  Thus, at least one overlap fails to resolve and the minimal generating set has at most 2 relations of degree three.\\
		
		In all cases where the  Gr{\"o}bner basis has at least 3 degree three relations, we find that the minimal generating set has at most 2 degree three relations, so there is no AS-Ore extension with relation type $(2,2,2,2,2,3,3,3)$.			
	\end{proof}
\end{theorem}

There are algebras with relation types with 0, 1, and 2 degree three relations.  We begin by considering what relation types can be realized by enveloping algebras.

\begin{theorem}\label{EnvelopingExamples}
	There are enveloping algebras with degree type (1,1,1,1,2) and relation type (2,2,2,2,2) and (2,2,2,2,2,3,3), but no enveloping algebra with relation type (2,2,2,2,2,3).
	\begin{proof}
		An enveloping algebra with such a degree type can be taken to have $x_5>x_4>x_3>x_2>x_1$ with $deg(x_1)$=2 and is then defined by the relations
		\begin{align*}
		r_ {21} : x_ 2 x_ 1 &= x_ 1 x_ 2 \\
		r_ {31} : x_ 3 x_ 1 &= x_ 1 x_ 3 \\
		r_ {41} : x_ 4 x_ 1 &= x_ 1 x_ 4 \\
		r_ {51} : x_ 5 x_ 1 &= x_ 1 x_ 5\\
		r_ {32} : x_ 3 x_ 2 &= b_ 1 x_ 1 + x_ 2 x_ 3 \\
		r_ {43} : x_ 4 x_ 3 &= d_ 1 x_ 1 + x_ 3 x_ 4 \\
		r_ {42} : x_ 4 x_ 2 &= e_ 1 x_ 1 + x_ 2 x_ 4\\
		r_ {54} : x_ 5 x_ 4 &= g_ 1 x_ 1 + x_ 4 x_ 5 \\
		r_ {53} : x_ 5 x_ 3 &= h_ 1 x_ 1 + x_ 3 x_ 5 \\
		r_ {52} : x_ 5 x_ 2 &= i_ 1 x_ 1 + x_ 2 x_ 5.
		\end{align*}
		
		These relations are homogeneous, overlaps resolve (see \cite[Section 4]{SE}), $\sigma_j$ is the identity for all $j\geq1$, $\delta_j$ is linear for all $j\geq1$ and this is generated in degree one if at least 1 of $b_1$, $d_1$, $e_1$, $g_1$, $h_1$, or $i_1$ is nonzero.  In this case, by \Cref{PresentationEnv}, this is an enveloping algebra and it remains to find its relation type.
		
		By the symmetry of the relations, we may assume without loss of generality that $b_1$ is nonzero and write $\displaystyle x_1=\frac{-x_2x_3+x_3x_2}{b_1}$. We can now view the algebra as $\widetilde{A}$, something generated in degree one, so that the set of LTs of degree two relations is $\{x_4x_3, x_4x_2, x_5x_4, x_5x_3, x_5x_2\}$.  Given this set of degree two LTs, we see from the analysis preceding \Cref{No(3333)} that the \basis $ $ of $\widetilde{A}$ has 2 degree three leading terms.  Further, the degree three overlaps that come from this list of degree two LTs are $x_5x_4x_3$ and $x_5x_4x_2$.  We calculate: 
		\begin{align*}
		x_5(x_4x_3)&-(x_5x_4)x_3=(\frac{g_1}{b_1}-\frac{e_1h_1}{b_1^2}+\frac{d_1i_1}{b_1^2})x_2x_3x_3 \\
		&+(\frac{-2g_1}{b_1}+\frac{2e_1h_1}{b_1^2}-\frac{2d_1i_1}{b_1^2})x_3x_2x_3+(\frac{g_1}{b_1}-\frac{e_1h_1}{b_1^2}+\frac{d_1i_1}{b_1^2})x_3x_3x_2\text{, and}\\
		x_5(x_4x_2)&-(x_5x_4)x_2=(\frac{-g_1}{b_1}+\frac{e_1h_1}{b_1^2}-\frac{d_1i_1}{b_1^2})x_2x_2x_3 \\
		&+(\frac{2g_1}{b_1}-\frac{2e_1h_1}{b_1^2}+\frac{2d_1i_1}{b_1^2})x_2x_3x_2+(\frac{-g_1}{b_1}+\frac{e_1h_1}{b_1^2}-\frac{d_1i_1}{b_1^2})x_3x_2x_2.
		\end{align*} 
		 As usual, the details for these calculations are omitted here but included in the code posted online \cite[Section 4]{SE}.  If $g_1=\displaystyle \frac{e_1h_1}{b_1}-\frac{d_1i_1}{b_1}$ then both of these overlaps resolve, the degree three relations in the \basis $ $ are independent of overlaps, and the relation type is (2,2,2,2,2,3,3).  Otherwise, these are two $K$-independent overlaps that fail to resolve and the minimal relation type is (2,2,2,2,2).  Since any enveloping algebra must have one of these two relation types, the theorem is proven.		 
	\end{proof}	
\end{theorem}

Although there is no enveloping algebra with relation type (2,2,2,2,3), it is possible to construct an AS-Ore extension with this type.  Our process for doing this is similar to that used to find the AS-Ore extension of degree type (1,1,2,3,5) in \Cref{Ore11235}. We use \Cref{Presentation 1} to write general relations for the extension and a mathematical program to evaluate the possible values of coefficients that make it so that those relations satisfy the diamond condition.  This problem is generally too complex for the computer to handle, so we can set some coefficients equal to zero to simplify the process.  In this case, we also have to determine how many degree three relations the \basis $ $ of the algebra generated in degree one has, as well as whether these relations are part of the minimal generating set as opposed to consequences of overlaps that do not resolve.

\begin{theorem}
	There is an AS-Ore extension with degree type (1,1,1,1,2) and relation type (2,2,2,2,2,3).
	\begin{proof}
		Consider the algebra defined by the relations
		\begin{align*}
		r_ {21} : x_ 2 x_ 1 & = x_ 1 x_ 2 \\
		r_ {32} : x_ 3 x_ 2 & = x_ 1 + x_ 2 x_ 2 - x_ 2 x_ 3 \\
		r_ {31} : x_ 3 x_ 1 & = x_ 1 x_ 3 \\
		r_ {43} : x_ 4 x_ 3 & = x_ 2 x_ 2 - x_ 3 x_ 4 \\
		r_ {42} : x_ 4 x_ 2 & = -x_ 2 x_ 4 \\
		r_ {41} : x_ 4 x_ 1 & = x_ 1 x_ 4 \\
		r_ {54} : x_ 5 x_ 4 & = -x_ 3 x_ 3 - x_ 4 x_ 5 \\
		r_ {53} : x_ 5 x_ 3 & = x_ 3 x_ 3 - x_ 3 x_ 5 \\
		r_ {52} : x_ 5 x_ 2 & = -x_ 2 x_ 5 + x_ 3 x_ 3 \\
		r_ {51} : x_ 5 x_ 1 & = x_ 1 x_ 5. \\
		\end{align*}
		Taking $x_5>\cdots>x_1$, the leading terms are as presented and all overlaps resolve \cite[Section 5]{SE}.  Taking $deg(x_1)=2$ and $deg(x_i)=1$, $2\leq i\leq5$, the relations are homogeneous.  Thus by \Cref{Presentation2}, this is an Ore extension.  It is also generated in degree one since $r_{32}$ can be used to express $x_1$ in terms of the degree one generators.  Finally, $\sigma_j(x_i)=\pm 1$ for all $1\leq i<j\leq5$ so these maps are injective and this is an AS-Ore extension.
		
		We can solve $r_{32}$ for $x_1$ and then view the algebra as $\widetilde{A}$, generated in degree one.  From the analysis preceding \Cref{No(3333)}, the \basis $ $ of $\widetilde{A}$ has 2 degree three relations and it remains to show that exactly one of these is a consequence of an overlap that fails to resolve.  We compute
		\[	x_5(x_4x_3)-(x_5x_4)x_3=x_2x_2x_3-x_2x_3x_3-x_3x_2x_2+x_3x_3x_2\]
		so this overlap never resolves.  Reducing the remaining overlap modulo this additional relation,
		\[x_5(x_4x_2)-(x_5x_4)x_2=0.\] 
		Thus, 1 of the degree three relations in the \basis $ $ is minimal and the relation type of this algebra is (2,2,2,2,2,3).
	\end{proof}
\end{theorem}

\section{A classification of relation types for AS-Ore extensions with 3 degree one generators}

We now begin the process of classifying the relation types of AS-Ore extensions with 3 degree one generators.  Again, by the symmetry of the free resolution, this also provides us with all the information we need to classify all possible resolution types of such algebras.  We recall that, by \Cref{DegreeType}, there are two possible degree types for algebras with 3 generators: (1,1,1,2,2) and (1,1,1,2,3).

 \begin{theorem} \label{(11122)}
An AS-Ore extension with degree type (1,1,1,2,2) has relation type (2,3,3,3,3,3).

\begin{proof}
By \Cref{Resolution}, a dimension 5 AS-Ore extension generated in degree one by 3 generators has free resolution
\[0\rightarrow A(-l)\rightarrow A(-l+1)^b\rightarrow \bigoplus \limits_{i=1}^{n} A(-l+a_i)\rightarrow \bigoplus \limits_{i=1}^{n}  A(-a_i)\rightarrow A(-1)^b\rightarrow A\rightarrow K\rightarrow 0\]
where b=3 since there are 3 degree one generators, $n$ represents the number of relations in the minimal generating set and $a_i$ represents the degree of the $i$th relation.  This algebra has Hilbert series 
\[h_A(t)=\frac{1}{q(t)}\text{  where } q(t)=1-3t +\sum \limits_{i=1}^{n} t^{a_i}-\sum \limits_{i=1}^{n} t^{l-a_i}+3t^{l-1}-t^l.\]

\iffalse By [stephenson zhang growth of graded noetherian rings], the Gelfand-Kirillow dimension is the order of the pole of $q(t)$ at $t=1$ so we get that
\begin{align}
\nonumber 0&=q'(1)\\
\nonumber &= -b+\sum\limits_{i=1}^n a_i-\sum\limits_{i=1}^n l-a_i+b(l-1)-l\text{ or} \nonumber\\
\label{firstderivative} \sum\limits_{i=1}^n a_i&=\frac{2b+(n-b+1)l}{2}, \\
\nonumber 0&=q''(0) \\
\nonumber &=\sum\limits_{i=1}^na_i(a_i-1)-\sum\limits_{i=1}^n(l-a_i)(l-a_i-1)+b(l-1)(l-2)-l(l-1)\\
\nonumber &= \sum\limits_{i=1}^na_i^2-a_i-\sum\limits_{i=1}^n a_i+a_i^2-l-2a_il+l^2+b(l-1)(l-2)-l(l-1)\\
\nonumber &=-nl^2+nl+(2l-2)\frac{2b(n-b+1)}{2}+b(l-1)(l-2)-l(l-1) \text{ by \ref{firstderivative}}\\
\nonumber &=0,\\
\nonumber 0&=q^{(3)}(1)\\
\nonumber &=\sum\limits_{i=1}^na_i(a_i-1)(a-i-2)-\sum\limits_{i=1}^n(l-a_i)(l-a_i-1)(l-a1-2)+b(l-1)(l-2)(l-3)-l(l-1)(l-2)\\
\nonumber &=\sum\limits_{i=1}^n(a_i^3-3a_i^2+2a_i)-\sum\limits_{i=1}^n(-2a_i-3a_i^2-a_i^3+2l+6a_il+3a_i^2l-3l^2-3a_il^2+l^3)\\
\nonumber & \;\;\; +b(l-1)(l-2)(l-3)-l(l-1)(l-2)
\end{align}
\fi

On the other hand, an AS-Ore extension with degree type (1,1,1,2,2) has Hilbert series
\[ h_A(t)=\frac{1}{(1-t)^3(1-t^2)^2} \text{, so } q(t)=1 - 3 t + t^2 + 5 t^3 - 5 t^4 - t^5 + 3 t^6 - t^7.\]
Assume $a_1\leq \cdots \leq a_n$.  Then $l=7$, $a_1=2$, $a_i=3$ for $2\leq i\leq 6$ (there are 5 degree three relations), and if there are any other minimal relations, they must cancel in the expression $\sum \limits_{i=1}^{n} t^{a_i}-\sum \limits_{i=1}^{n} t^{l-a_i}$ since they do not appear in the second equation for $q(t)$.  This would mean either that $a_i=l-a_i$ (which is impossible since $l$ is odd) or that there are at least two additional relations, $a_i$ and $a_j$ with $a_i+a_j=l$.  Although we are interested in the minimal generating set of $\A$, the algebra generated in degree one, we note that any minimal relations of $\A$ must come from the original relations of the algebra viewed as an Ore extension, $A$.  The possible relations for $A$ are described in the presentation of an Ore extension given by \Cref{Presentation 1} and an Ore extension with degree type (1,1,1,2,2) can only have relations of degrees two, three, or four.  Thus, if there are additional minimal relations that cancel in the Hilbert series, they must be of degree three and four (since something of degree two or lower could only cancel if there were also a relation of degree five or higher, and we know that an Ore extension with this degree type has no such relations which are minimal).  

If we label the degree one generators with the order $x_3<x_4<x_5$, the algebra, when viewed as $\A$, will have 1 degree two relation with leading term $x_5x_4,$ $x_5x_3,$ or $x_4x_3$. (The remaining 2 degree two relations will be used to write the $x_1$ and $x_2$ in terms of the degree one generators.)  If $\A_2$ denotes the monomial algebra generated in degree one that has one of $x_5x_4,$ $x_5x_3,$ or $x_4x_3$ as a relation and $\A$ denotes the AS-Ore extension with degree type (1,1,1,2,2) viewed as an algebra generated in degree one, then $h_{\A_2}(t)-h_{\A}(t)=5t^3+O(t^4)$ \cite[Section 6]{SE}, so there can only be 5 degree three relations in the \basis.  Thus, the only relation type for an AS-Ore extension with degree type (1,1,1,2,2) is (2,3,3,3,3,3).
\end{proof}
\end{theorem}

\begin{theorem} \label{Env(11122)}
There is an enveloping algebra with degree type (1,1,1,2,2) and relation type (2,3,3,3,3,3).

\begin{proof}
Consider the algebra defined by relations
\begin{align*}
r_ {21} : x_ 2 x_ 1 &= x_ 1 x_ 2\\
r_ {32} : x_ 3 x_ 2 &= x_ 2 x_ 3\\
r_ {31} : x_ 3 x_ 1 &= x_ 1 x_ 3\\
r_ {42} : x_ 4 x_ 2 &= x_ 2 x_ 4\\
r_ {41} : x_ 4 x_ 1 &= x_ 1 x_ 4\\
r_ {52} : x_ 5 x_ 2 &= x_ 2 x_ 5\\
r_ {51} : x_ 5 x_ 1 &= x_ 1 x_ 5\\
r_ {43} : x_ 4 x_ 3 &= x_ 1 + x_ 3 x_ 4\\
r_ {54} : x_ 5 x_ 4 &= x_ 2 + x_ 4 x_ 5\\
r_ {53} : x_ 5 x_ 3 &= x_ 3 x_ 5.\\
\end{align*}
Assigning $(x_5,x_4,x_3,x_2,x_1)$ degrees $(1,1,1,2,2)$, these relations are homogeneous and it can be verified by hand or computer that all overlaps resolve, so this is an Ore extension by \Cref{Presentation2}. Additionally, for all $1\leq i<j\leq 5$, $\sigma_j(x_i)$ is the identity and $\delta_j(x_i)$ is linear, so this is an enveloping algebra by \Cref{PresentationEnv}.  It is also generated in degree one.  Another quick check in Mathematica shows that the relation type is $(2,3,3,3,3,3)$ \cite[Section 7]{SE}, although the analysis from \Cref{(11122)} already indicates that this has to be the case since this is the only possible relation type for an AS-Ore extension with this degree type.
\end{proof}
\end{theorem}

\begin{theorem} \label{(11123)}
An AS-Ore extension with degree type (1,1,1,2,3) has relation type (2,2,3) or (2,2,3,4).
\begin{proof}
	Following the logic of \Cref{(11122)}, we know from the free resolution of the algebra that
	\[h_A(t)=\frac{1}{q(t)}\text{  where } q(t)=1-3t +\sum \limits_{i=1}^{n} t^{a_i}-\sum \limits_{i=1}^{n} t^{l-a_i}+3t^{l-1}-t^l.\]
	On the other hand, an AS-Ore extension with degree type (1,1,1,2,3) has Hilbert series
	\[ h_A(t)=\frac{1}{(1-t)^3(1-t^2)(1-t^3)} \text{, so } q(t)=1 - 3 t + 2 t^2 + t^3 - t^5 - 2 t^6 + 3 t^7 - t^8.\]
	So $l=8$, $a_1=a_2=2$, $a_3=3$, and if there are any other minimal relations of $\A$, generated in degree one, they must cancel.  In this case we would have $a_i=l-a_j$ (where it is possible that $i=j$).  This means that there may possibly be a pair of minimal relations of degree three and five (there cannot be more than one such pair since there is only 1 degree five relation in the Ore extension) and there may possibly be up to 3 relations of degree four (since the Ore extension $A$ has 3 relations of degree four).  There cannot be any additional minimal relations of degree two since they would need to cancel with something of degree six and there are no minimal relations of degree six for an AS-Ore extension with this degree type.   To conclude that the relation type is either (2,2,3) or (2,2,3,4), it remains to show that there are not 2 independent relations of degree three (and so no additional pair of relations of degree three and five) and that there is at most 1 minimal relation of degree four.
	
	We can write relations based on the fact that this is an Ore extension.  Without loss of generality, label the degree one generators so that $x_3<x_4<x_5$.  We will let $x_2$ be the degree two variable and $x_1$ will be the degree three variable, but we make no claim about when these variables are adjoined.  (Thus possible orders include $x_1<x_2<x_3<x_4<x_5$ and $x_3<x_2<x_4<x_1<x_5$, amongst many others.)  We do note, however, that for $A$ to be generated in degree one, $x_1$ must not be adjoined last (since it must be adjoined by the time it appears in a relation with leading term of the form $x_jx_i$ and hence must be adjoined before $x_j$), and (similarly) $x_2$ must not be adjoined last.  Then by the choice of ordering of our degree one variables, $x_4$ and $x_3$ are both adjoined before $x_5$, so $x_5$ is always added last in the iterated Ore extension.
	
	Using the ordering of the degree one variables, we can write the degree two relations.  We note that the only thing that could change in the relations below that depends on the order in which the variables is adjoined is that $d_1$ must be zero if $x_2$ is added after $x_4$, since $r_{43}$ should only involve variables that have been added by the time $x_4$ is adjoined.  The degree two relations are:
	\begin{align*}
	r_{43}: x_4x_3 & = d_ 1 x_ 2 +d_ 2 x_ 3 x_ 3 + d_ 3 x_ 3 x_ 4\\
	r_{53}: x_5x_3 & = h_ 1 x_ 2 + h_ 2 x_ 3 x_ 3 + h_ 3 x_ 3 x_ 4+    h_ 4 x_ 3 x_ 5 + h_ 5 x_ 4 x_ 4+ h_ 6 x_ 4 x_ 5 \\
	r_{54}: x_5x_4 & = g_ 1 x_ 2 + g_ 2 x_ 3 x_ 3 + g_ 3 x_ 3 x_ 4+    g_ 4 x_ 3 x_ 5 + g_ 5 x_ 4 x_ 4 + g_ 6 x_ 4 x_ 5.
	\end{align*}
	We wish to repeat this process for the degree three relations. The list of possible degree three monomials that can appear on the right side of a relation, given the ordering $x_3<x_4<x_5$, $x_1<x_5$, $x_2<x_5$, is:
	\[\{x_1,x_2x_3,x_3x_2,x_2x_4,x_4x_2,x_2x_5,x_3x_3x_3,x_3x_3x_4,x_3x_3x_5,x_3x_4x_4,x_3x_4x_5, x_4x_4x_4\}.\]
	Since the ordering of the variables is not known, the leading terms of the degree 3 relations can vary, as can the other allowable monomials. 
	We will write fully general versions of the relations for each possible leading term.  It is appropriate to use $r_{32a}$ below when $x_3>x_2$ and $r_{32b}$ when $x_2>x_3$.  Similarly, $r_{42a}$ applies when $x_4>x_2$, and $r_{42b}$ when $x_2>x_4$.
		
	 From \Cref{Presentation 1}, $x_j$ must occur only to the first power in the relation with leading term $x_jx_i$ and $x_k$ should not appear for any $x_k>x_j$.
	\iffalse
	For example in $r_{32a}$ below, having a leading term of $x_3x_2$ implies that $x_4>x_3>x_2$ so no monomial containing $x_4$ will appear in the relation.  Also note that for $r_{42b}$, $x_2>x_4>x_3$ implies that $x_2x_3$ will not appear since it can be rewritten in terms of smaller monomials.
	\fi 
	For example in $r_{32b}$ below, the leading term implies that $x_2>x_3$.  If $x_2>x_4>x_3$ then the monomial $x_4x_2$ may appear.  Otherwise, $x_4>x_2>x_3$ and the monomial $x_2x_4$ does not appear since $x_4$ has not yet been adjoined in the Ore extension.  The most general possible degree three relations are:
	\begin{align*}	
	r_{32}a:&&-b_0x_3x_2&=b_1x_1+b_2x_2x_3\\
	r_{32}b:&&-b_2x_2x_3&=b_1x_1+b_0x_3x_2+b_3x_4x_2+b_4x_3x_3x_3+b_5x_3x_3x_4+b_6x_3x_4x_4\\
	&&&+b_7x_4x_4x_4\\
	\end{align*}
	\begin{align*}
	r_{42}a:&&-e_0x_4x_2&=e_1x_1+e_2x_2x_3+e_3x_2x_4+e_4x_3x_2+e_5x_3x_3x_3+e_6x_3x_3x_4\\
	r_{42}b:&&-e_3x_2x_4&=e_1x_1+e_0x_4x_2+e_4x_3x_2+e_5x_3x_3x_3+e_6x_3x_3x_4+e_7x_3x_4x_4\\
	&&&+e_8x_4x_4x_4\\
	r_{52}:&&x_5x_2&=i_ 1 x_ 1 + i_ 2 x_ 2 x_ 3+ i_ 3 x_ 2 x_ 4 +i_4x_3x_2+i_5x_4x_2+    i_6 x_ 2 x_ 5\\
	&&& +   i_7 x_ 3 x_ 3 x_ 3 +    i_ 8 x_ 3 x_ 3 x_ 4+i_ 9 x_ 3 x_ 3 x_ 5 + i_ {10} x_ 3 x_ 4 x_ 4 +    i_{11} x_ 3 x_ 4 x_ 5\\
	&&& + i_{ 12} x_ 4 x_ 4 x_ 4+i_{13}x_4x_4x_5.
	\end{align*}
	Note that not all of these coefficients can be nonzero.  For example, if $x_2<x_4$ then $i_{5}=0$ while if $x_4<x_2$, $i_3=0$.  These do, however, capture all possible terms that could occur in the Ore relations.
	
	\iffalse
	In the work that follows, We will use the fact that several of the coefficients above cannot be zero.  These coefficients include $i_0$ (the leading term of $r_{52}$, regardless of the ordering of the variables), $d_3$ (by the injectivity of $\sigma_3$),  $e_3$ (for $r_{42}a$, $\sigma_4(x_2)=e_3x_2+e_6x_3x_3$, $\sigma_4(x_3)=d_3x_3$ by $r_{43}$ and so $\sigma_4(x_2-\displaystyle\frac{e_6}{d_3^2}x_3x_3)=e_3x_2+e_6x_3x_3-e_6x_3x_3$ and so $e_3\neq0$ by the injectivity of $\sigma_4$.  For $r_{42}b$, $e_3$ is the leading term.), and $b_2$ (for $r_{32}a$, $\sigma_3(x_2)=b_2x_2x_3$ and $b_2\neq0$ by the injectivity of $\sigma_3$.  Otherwise, $b_2$ is the leading term of $r_{32}b$.)
	\fi
	
	The goal now is to view this as $\A$, generated in degree one, and to try to identify the possible degrees of relations that can occur.  The proof depends on which relations are used to solve for $x_2$.  Without loss of generality, we will solve for $x_1$ and $x_2$ using the relation with smallest leading term so that the leading terms of the remaining degree two and three relations are easily identifiable. 
	% By using the automorphisms $x_2\to\frac{1}{k_2}x_2$ and $x_1\to\frac{1}{k_1}x_1$, $k_1,k_2\in K^\times$, we may also assume without loss of generality that the coefficients in the relations used to solve for $x_2$ and $x_1$ are 1.
	\\
	
	\textbf{Case 1:} $x_2$ comes from $r_{43}$.\\
	In this case, the degree two leading terms are $x_5x_4$ and $x_5x_3$.  If $\A_2$ is the monomial algebra with these leading terms then $h_{\A_2}-h_{\A}=t^3+O(t^4)$ so there is 1 degree three relation in the \basis $ $ of $\A$ and so at most (and exactly) 1 degree three relation in the minimal generating set, as we wished to show.  It remains to show that there is no more than 1 degree four relation in the minimal generating set.
	
	Since $x_2$ comes from $r_{43}$, $x_2$ certainly appears in the relation $r_{43}$ and so must have been adjoined before $x_4$.  This means that $x_2<x_4$ in the order which in turn implies that we should use the relation $r_{42}a$, that $e_0$ is not zero since it is the leading term of the expression, and that $b_3=b_5=b_6=b_7=0$.  We can begin to view the algebra as generated in degree one by solving $r_{43}$ to get $\displaystyle x_2=\frac{1}{d_1}(x_4x_3-d_3x_3x_4-d_2x_3x_3)$.
	We can substitute this into the other relations and move all terms to the right side of the equation:
	\begin{align*}
	r_{32a}:0&=b_0 x_3x_2+b_1x_1+b_2x_2x_3\\
	&=b_0 x_3\frac{1}{d_1}(x_4x_3-d_3x_3x_4-d_2x_3x_3)+b_1x_1\\
	& \phantom{we}+b_2\frac{1}{d_1}(x_4x_3-d_3x_3x_4-d_2x_3x_3)x_3\\
	&=b_1 x_1+\frac{b_2}{d_1}x_4x_3x_3+\text{[smaller terms]}\\
		r_{32b}:0&=b_2x_2x_3+b_1x_1+b_0x_3x_2+0x_4x_2+b_4x_3x_3x_3+0x_3x_3x_4+0x_3x_4x_4\\
	&\phantom{wo}+0x_4x_4x_4\\
	&=b_2\frac{1}{d_1}(x_4x_3-d_3x_3x_4-d_2x_3x_3)x_3+b_1x_1+b_4x_3x_3x_3\\
	&=b_1x_1+\frac{b_2}{d_1}x_4x_3x_3+\text{[smaller terms]}.
	\end{align*}

	 We note that the largest term appearing in the relation for $r_{32}$ is now independent of which version of the relation we use and, repeating the substitution for $x_2$ in $r_{42a}$, we can write
	\begin{align*}
	r_{32}:0&=b_1x_1+\frac{b_2}{d_1}x_4x_3x_3+\text{[smaller terms]}\\
	r_{42}:0&=e_1x_1+\frac{e_0}{d_1}x_4x_4x_3+\text{[smaller terms]}.
	\end{align*}
	From the original relations, we observe that $b_2$ is not zero: either $x_2>x_3$ in the ordering and $b_2$ is the LT of $r_{32}$ or $x_3>x_2$ and $\sigma_3(x_2)=b_2x_2$ and so $b_2\neq0$ by the injectivity of $\sigma_3$.
	
	If $b_1$ is not zero, we may solve $r_{32}$ for $x_1$ and, even after substituting this value into $r_{42}$, the leading term of $r_{42}$ in $\A$ will be	$x_4x_4x_3$.  If $b_1$ is zero then $x_4x_3x_3$ will be a leading term in $\A$.  If $\A_3$ is the monomial algebra with LTs $\{x_5x_4, x_5x_3, x_4x_3x_3\}$ or $\{x_5x_4, x_5x_3, x_4x_4x_3\}$ then $h_{\A_3}-h_{\A}=t^4+O(t^5)$.  So the \basis $ $ of $\A$ has 1 degree four relation and so at there is at most 1 degree four relation in the minimal generating set, as we wished to show.\\
	
	\textbf{Case 2:} $x_2$ comes from $r_{54}$.\\
	In this case, the degree two leading terms are $x_5x_3$ and $x_4x_3$, $d_1=g_1=0$ (or else $x_2$ would come from $r_{43}$ or $r_{53}$), and $h_{\A_2}-h_{A}=t^3+O(t^4)$.  This means there is at most 1 degree three relation in the minimal generating set and it remains to show that there is also at most 1 degree four relation.  From $r_{54}$, we can write $x_2=\displaystyle \frac{1}{g_1}(x_5x_4+\text{ [smaller terms]})$.  As in the first case, we can substitute this value into the degree three terms to get that the highest terms in the degree three relations of interest are known, independent of which version of $r_{42}$ we use:
	\begin{align*}
	r_{42}:0&=e_1x_1+ \frac{e_3}{g_1}x_5x_4x_4+\text{[smaller terms]}\\
	r_{52}:0&=i_1x_1-\frac{1}{g_1}x_5x_5x_4+\text{[smaller terms]}.
	\end{align*}
	Note that $e_3$ is not zero: either $x_4<x_2$ and it $e_3$ the leading coefficient of $r_{42b}$, or $x_2<x_4$ and from $r_{42a}$,
	$\sigma_4(x_2)=\displaystyle -\frac{e_3}{e_0}x_2-\frac{e_6}{e_0}x_3x_3$.  In this case, since $\sigma_4(x_3)=d_3x_3$, $\sigma_4(x_2+\displaystyle \frac{e_6}{e_0d_{3}^2}x_3x_3)=-\frac{e_3}{e_0}x_2-\frac{e_6}{e_0}x_3x_3+\frac{e_6}{e_0}x_3x_3$. By the injectivity of $\sigma_4$, $e_3\neq0$.  
	
	If $e_1$ is not zero, we may solve $r_{42}$ for $x_1$ and, even after substituting this value into $r_{52}$, the LT of $r_{52}$ will be $x_5x_5x_4$.  If $e_1$ is zero then $x_5x_4x_4$ will be a leading term in $\A$.  In either case, we calculate that $h_{\A_3}-h_{\A}=t^4+O(t^5)$, which means that there is at most 1 degree four relation in the minimal generating set, as we wished to prove.\\
	
	\textbf{Case 3:} $x_2$ comes from $r_{53}$.\\
	In this case, the degree two LTs are $x_5x_4$ and $x_4x_3$, $d_1=0$, and $h_{\A_2}-h_{\A}=2t^3-O(t^4)$.  We can solve $r_{53}$ to get that $x_2=\displaystyle \frac{1}{h_1}(x_5x_3-h_6x_4x_5+\text{[smaller terms]})$.  We can rewrite the remaining degree two relations after substituting the value of $x_2$ into the equations:
	\begin{align*}
	r_{43}:x_4x_3&=0x_2+d_3x_3x_4+\text{[smaller terms]}\\
	r_{54}:x_5x_4&=g_1\frac{1}{h_1}(x_5x_3+\text{[smaller terms]})+g_6x_4x_5+\text{[smaller terms]}.
	\end{align*}
	We can then compute the degree three overlap in $\A$.
	\begin{align*}
	(x_5x_4)x_3-x_5(x_4x_3)&=\frac{g_1}{h_1}(x_5x_3-h_6x_4x_5+\text{[small]})x_3-x_5(d_3x_3x_4+d_2x_3x_3)\\
	&=-d_3x_5x_3x_4+(\frac{g_1}{h_1}-d_2)x_5x_3x_3+\text{[small]}.
	\end{align*}
	Since $d_3$ is not zero by the injectivity of $\sigma_3$, this overlap does not resolve and we know that there is a relation, $x_5x_3x_4=\displaystyle \frac{g_1-d_2h_1}{h_1}x_5x_3x_3+\text{[smaller terms]}$ in the \basis $ $ of $\A$.  We conclude that at most 1 of the 2 degree three relations in the \basis $ $ can be minimal and it remains to show that there is at most 1 minimal relation of degree four.
	
	Again substituting the value of $x_2$, we may re-examine the degree three relations and note that, as with the first case, the LT of $r_{32}$ is independent of which version we use:
	\begin{align*}
	r_{32}:0&=b_1x_1+\frac{b_2}{h_1}x_5x_3x_3+\text{[smaller terms]}\\
	r_{52}:0&=i_1x_1-\frac{1}{h_1}x_5x_5x_3+\text{[smaller terms]}.
	\end{align*}
	By the same analysis as in case 1, $b_2\neq0$.  If $b_1=0$ then the leading term of $r_{32}$ in $\A$ is $x_5x_3x_3$.  If $\A_3$ is the monomial algebra with LTs $\{x_5x_4, x_4x_3, x_5x_3x_3,x_5x_3x_4\}$ then $h_{\A_3}-h_{\A}=t^4+O(t^5)$ and there is at most 1 degree 4 relation in the minimal generating set as desired.
	
	If $b_1$ is not zero then $r_{32}$ may be solved for $x_1$ and, even after substituting the value of $x_1$ into the relation, the leading term of $r_{52}$ in $\A$ is $x_5x_5x_3$.  Thus, the 2 relations of degree three in the \basis $ $ have LTs $x_5x_3x_4$ and $x_5x_5x_3$ and in this case, $h_{\A_3}-h_{\A}=2t^4+O(t^5)$.  We can also compute the degree four overlap:

	\begin{equation*} 
	\begin{split}
	(x_5x_3x_4)x_3-x_5x_3(x_4x_3)&=(\displaystyle\frac{g_1-d_2h_1}{d_3h_1}x_5x_3x_3+\text{[small]})x_3 \\
		&\hspace{-2 cm} -x_5x_3(d_3x_3x_4+\text{[small]})	\\	
	& =-d_3x_5x_3x_3x_4+\text{[small]}.
	\end{split}
	\end{equation*}
	
	We note that $x_5x_3x_3x_4$ cannot be reduced in $\A_3$ and $d_3$ is not zero so this overlap fails to resolve.  In total, we have found that the \basis $ $ of $\A$ has 2 degree four relations, at least one of which is not minimal. \\ 
	
	Thus, in all cases we have shown that $\A$ has at most 1 degree three and at most 1 degree four relation in the minimal generating set, which means that the relation type must be either (2,2,3), or (2,2,3,4).	
	\end{proof}

\end{theorem}

\begin{theorem} \label{Env(11123)}
	There is an enveloping algebra with degree type (1,1,1,2,3) and relation type (2,2,3,4), but not one with relation type (2,2,3).

\begin{proof}
	An enveloping algebra can be taken to have $x_5>x_4>x_3>x_2>x_1$ with $deg(x_2)=2$ and $deg(x_1)=3$ and is then defined by the relations
	\begin{align*}
	r_ {21} : x_ 2 x_ 1 &= x_ 1 x_ 2 \\
	r_ {31} : x_ 3 x_ 1 &= x_ 1 x_ 3 \\
	r_ {41} : x_ 4 x_ 1 &= x_ 1 x_ 4 \\
	r_ {51} : x_ 5 x_ 1 &= x_ 1 x_ 5 \\
	r_ {32} : x_ 3 x_ 2 &= b_1 x_1+x_ 2 x_ 3 \\
	r_ {42} : x_ 4 x_ 2 &= e_1 x_1+x_ 2 x_ 4 \\
	r_ {52} : x_ 5 x_ 2 &= i_1 x_ 1 + x_ 2 x_ 5 \\
	r_ {43} : x_ 4 x_ 3 &= d_1 x_2+x_ 3 x_ 4 \\
	r_ {54} : x_ 5 x_ 4 &= g_1 x_2+ x_ 4 x_ 5 \\
	r_ {53} : x_ 5 x_ 3 &= h_1 x_2+ x_ 3 x_ 5. \\
	\end{align*}
	
By construction we have that for all $1\leq i<j \leq 5$, $\sigma_j(x_i)$ is the identity and $\delta_j(x_i)$ is linear for all $i$.  All overlaps resolve, except for $x_5(x_4x_3)-(x_5x_4)x_3=(b_1 g_1 - e_1 h_1 + d_1 i_1) x_1$, so we will have to choose values of coefficients which make this expression zero.  Additionally, to be generated in degree one, we must have that at least 1 of $b_1$, $e_1$, $i_1$ and at least 1 of $d_1$, $g_1$, $h_1$ is nonzero.  If this happens, then by \Cref{PresentationEnv}, this is an enveloping algebra.

   We can now solve for $x_1$ and $x_2$ to view this as $\A$ and analyze the possible degrees of minimal relations. As the process for these computations is quite similar to that seen in previous theorems, we will omit most of the details.  Our goal is to show that the relation type is always (2,2,3,4). By the symmetry of the relations, we may assume that $b_1$ is the coefficient that is not zero.  \\

\textbf{Case 1:} $d_1$ nonzero.\\
From the overlap in $A$, $x_5(x_4x_3)-(x_5x_4)x_3=(g_1 b_1-e_1 h_1+i_1 d_1)x_1$, we conclude that $g_1=\displaystyle\frac{e_1 h_1-i_1 d_1}{b_1}$.  Solving $r_{43}$ and $r_{32}$ for $x_2$ and $x_1$ and substituting these into the remaining relations to view the algebra as generated in degree one, we find that there are degree two relations in $\A$ with LTs $x_5x_3$ and $x_5x_4$ and a degree three relation from $r_{42}$ with LT $x_4x_4x_3$.  Reduced modulo these relations, $r_{51}$ then has LT $x_4x_3x_3x_3$ and it remains to show that this is minimal in $\A$.  The only degree four overlap, given these LTs, is $(x_5x_4)x_4x_3-x_5(x_4x_4x_3)=0$.  Details for computations in this proof are available online \cite[Section 9]{SE}.  Since this overlap resolves, the degree four relation is independent.  By \Cref{(11123)}, this means that the relation type must be (2,2,3,4).  In particular, we have now shown that there is an enveloping algebra with this relation type.  It remains to show that (2,2,3) is never a possible relation type.\\

\textbf{Case 2:} $d_1=0$ and $h_1$ nonzero.\\
From the overlap $x_5(x_4x_3)-(x_5x_4)x_3=(-e_1 h_1+g_1 b_1)x_1$ we conclude that $e_1 = \displaystyle\frac{g_1 b_1}{h_1}$.  Solving $r_{53}$ and $r_{32}$ for $x_2$ and $x_1$ and substituting these into the remaining relations to view the algebra as generated in degree 1, we find that there are degree two relations in $\A$ with LTs $x_4x_3$ and $x_5x_4$, a degree three overlap that fails to resolve with LT $x_5x_3x_4$, and a degree three relation from $r_{52}$ with LT $x_5x_5x_3$.  The monomial algebra with these LTs has 2 degree four relations so it remains to show that exactly one such relation is the consequence of an overlap that does not resolve.  One such overlap is $x_5x_3(x_4x_3)-(x_5x_3x_4)x_3=\displaystyle \frac{g_1}{h_1}x_3x_5x_3x_3-x_4x_5x_3x_3-\frac{g_1}{h_1}x_5x_3x_3x_3+x_5x_3x_3x_4$ and so fails to resolve.  The remaining degree four overlap, when reduced modulo the degree two and three relations together with this new relation, is $x_5(x_5x_3)x_4-(x_5x_5x_3)x_4=0$.  Thus, there is 1 independent degree four relation and the relation type is (2,2,3,4).\\

\textbf{Case 3:} $d_1=0$, $h_1=0$, and $g_1$ nonzero.\\
Recall that, by the symmetry of the relations, we have assumed that $b_1\neq0$, and that $g_1\neq0$ if $d_1=0$, $h_1=0$, and $A$ is generated in degree one.  From the overlap in $A$, $x_5(x_4x_3)-(x_5x_4)x_3=b_1 g_1 x_1$, we conclude that there is no enveloping algebra with coefficients with these values since the overlaps of the original Ore relations fail to resolve.\\

Thus in all cases, the only possible relation type is (2,2,3,4) and so this is the only relation type of an enveloping algebra with variables of degrees (1,1,1,2,3).
\end{proof}
\end{theorem}

\begin{theorem}
	There is an AS-Ore extension with relation type (2,2,3).
	
\begin{proof}
	Consider the algebra defined by the relations
	\begin{align*}
	r_{21}: x_2x_1&=-x_1x_2\\
	r_{32}: x_3x_2&=x_2x_3\\
	r_{31}: x_3x_1&=x_1x_3\\
	r_{43}: x_4x_3&=x_3x_4\\	
	r_{42}: x_4x_2&=x_1+x_2x_4\\
	r_{54}: x_5x_4&=x_2+x_4x_5\\
	r_{53}:x_5x_3&=x_3x_5+x_4x_4\\
	r_{52}:x_5x_2&=x_2x_5\\
	r_{51}:x_5x_1&=x_1x_5.
	\end{align*}

Assigning $(x_5,x_4,x_3,x_2,x_1)$ degrees (1,1,1,2,3), these relations are homogeneous.  Using the order $x_5>x_4>x_3>x_2>x_1$ so that the relations are as presented, all overlaps resolve \cite[Section 10]{SE}.  Hence, this is an Ore extension by \Cref{Presentation2}.  It is also generated in degree one.  For all $1\leq i<j\leq 5$, $\sigma_j(x_i)=\pm 1$ so the $\sigma_j$ are automorphisms and this algebra is AS-Ore.
	
We can view this algebra as $\A$, generated in degree one, by solving $r_{54}$ and $r_{42}$ for $x_2$ and $x_1$ and plugging these values into the relations.  The remaining degree two relations in $\A$ have LTs $x_5x_3$ and $x_4x_3$ and there is a degree three relation with LT $x_5x_5x_4$.  Reduced modulo these relations, $r_{31}$ and $r_{51}$ become $0$ while $r_{41}=-x_4x_4x_4x_5+x_4x_4x_5x_4+x_4x_5x_4x_4-x_5x_4x_4x_4$.  So there is 1 degree four relation in the \basis $ $ (which we already knew from the Hilbert series analysis of \Cref{(11123)}) and it remains to show that this is not minimal.  We compute the overlap $x_5x_5(x_4x_3)-(x_5x_5x_4)x_3=-x_4x_4x_4x_5+x_4x_4x_5x_4+x_4x_5x_4x_4-x_5x_4x_4x_4$.  Thus, the degree four relation is a consequence of an overlap that fails to resolve and the relation type is (2,2,3).	
\end{proof}
\end{theorem}

\bibliographystyle{alphanum}
\bibliography{Bibliography}

\end{document}